\documentclass[a4paper]{article}

\usepackage[pdftex]{graphicx} 
\usepackage{amsfonts}
\usepackage{amsmath,amsthm}
\usepackage{latexsym}
\usepackage{array}
\usepackage{amssymb}
\usepackage{epstopdf} 
\usepackage{enumerate}
\usepackage[english]{babel}

\usepackage{mathtools}
\usepackage{graphics}
\usepackage{amssymb}
\usepackage{amsfonts}
\usepackage{amsmath,amsthm}

\usepackage{listings}
\usepackage{mathrsfs}
\usepackage{enumerate}
\usepackage{pict2e}
\usepackage{stmaryrd}
\usepackage{xcolor}

\usepackage[all,cmtip]{xy}

\usepackage{latexsym}
\usepackage{mathrsfs}
\usepackage[all,cmtip]{xy}

\usepackage{tikz}
\usepackage{amsthm}
\newtheorem*{remark}{Remark}

\usepackage{scalerel,stackengine}
\stackMath
\newcommand\Widehat[1]{%
\savestack{\tmpbox}{\stretchto{%
  \scaleto{%
    \scalerel*[\widthof{\ensuremath{#1}}]{\kern-.6pt\bigwedge\kern-.6pt}%
    {\rule[-\textheight/2]{1ex}{1.2\textheight}}
  }{\textheight}%
}{0.6ex}}%
\stackon[1pt]{#1}{\tmpbox}%
}
\newcommand\ie{\emph{i.e.}~}

\usepackage{algpseudocode,algorithm}
\algnewcommand{\LineComment}[1]{\State  \(\triangleright\) #1 \hfill~}

\usepackage{amsmath}

\usepackage{caption}
\usepackage{subcaption} 

\providecommand{\keywords}[1]{\textbf{\textit{Keywords---}} #1}
\usepackage{authblk} 

\title{Exponential methods for solving hyperbolic problems with application to collisionless kinetic equations}
\author[a]{Nicolas Crouseilles}
\author[b]{Lukas Einkemmer}
\author[c]{Josselin Massot}
\affil[a]{Univ Rennes, Inria Bretagne Atlantique (MINGuS) \& ENS Rennes}
\affil[b]{Department of Mathematics, University of Innsbruck, Austria}
\affil[c]{Univ Rennes \& Inria Bretagne Atlantique (MINGuS)}

\usepackage[pdftex,
            pdfauthor={Nicolas Crouseilles, Lukas Einkemmer, Josselin Massot},
            pdftitle={Exponential methods for solving hyperbolic problems with application to kinetic equations},
            pdfsubject={MSC:65},
            pdfkeywords={exponential integrators, Lawson schemes, kinetic equations, hyperbolic PDEs, numerical stability, drift-kinetic equations}]{hyperref}

\begin{document}
\maketitle

\abstract{The efficient numerical solution of many kinetic models in plasma physics is impeded by the stiffness of these systems. Exponential integrators are attractive in this context as they remove the CFL condition induced by the linear part of the system, which in practice is often the most stringent stability constraint. In the literature, these schemes have been found to perform well, {\it e.g.}, for drift-kinetic problems. Despite their overall efficiency and their many favorable properties, most of the commonly used exponential integrators behave rather erratically in terms of the allowed time step size in some situations. This severely limits their utility and robustness.

Our goal in this paper is to explain the observed behavior and suggest exponential methods that do not suffer from the stated deficiencies. To accomplish this we study the stability of exponential integrators for a linearized problem. This analysis shows that classic exponential integrators exhibit severe deficiencies in that regard. Based on the analysis conducted we propose to use Lawson methods, which can be shown not to suffer from the same stability issues. We confirm these results and demonstrate the efficiency of Lawson methods by performing numerical simulations for both the Vlasov--Poisson system and a drift-kinetic model of a ion temperature gradient instability.}

\ 


\section{Introduction}
\label{intro}

The goal of this work is to develop high order and efficient numerical methods for nonlinear \textcolor{black}{collisionless} kinetic models, such as the Vlasov-Poisson equations or drift-kinetic models. In most situations, the nonlinearity in the transport term originates from the coupling with a Poisson type problem that is used to compute the electric field.  

Historically, particle in cell methods have been extensively used to treat kinetic problems. In this approach, the unknown is sampled by discrete particles which are advanced in time using an ODE solver, whereas the electric field is computed on a spatial grid. For some problems, these methods can tackle high dimensional kinetic problems with a relatively low computational cost. However, they also suffer from numerical noise which pollutes the accuracy in low density regions of phase space. Moreover, as the number of particles is increased the error only decreases as the inverse of the square of the number of particles. Thus, convergence is slow. For a review of particle methods we refer the reader to \cite{verboncoeur2005}.

On the other hand, Eulerian methods ({\it e.g.}~finite volumes or finite differences), which directly discretize the phase space, are able to reach high order accuracy in time, space, and velocity.  However, in addition to the fact that they are costly, these methods usually suffer from stability constraints that force a relation between the time step and the phase space grid sizes, the so-called Courant--Friedrichs--Lewy (CFL) condition. Hence, a large number of time steps is required to reach the long times that are often required in plasma physics applications. 

To overcome this CFL condition, semi-Lagrangian methods have been developed during the last decades. These methods realize a compromise between the Lagrangian (\ie~particle in cell) and Eulerian approaches by exploiting the characteristics equations to overcome the CFL condition, while still performing computations on a grid in both space and velocity \cite{cheng1976,sonnendrucker1999,filbet2003}. This approach is usually combined with splitting methods to avoid a costly multidimensional interpolation step. This allows for a separate treatment of the terms in the equation and the corresponding characteristic curves can then, at least in some situations, be computed analytically. For purely hyperbolic problems, the setting we consider here, it is also possible to construct high order splitting schemes \cite{casas2017}.

Splitting results in a very accurate and efficient scheme for the Vlasov-Poisson equation. This is the case because the problem is only split into two parts, the characteristics of which can then be solved exactly in time.  However, this is not necessarily true for more complicated equations such as gyrokinetic or drift-kinetic models. Indeed, for the drift-kinetic model, a three terms splitting has to be performed so that a relatively large number of stages are required to reach high order accuracy in time.  In addition, some stages can not be solved exactly in time and thus require additional numerical work to approximate them.

In \cite{cep} an alternative approach based on exponential integrators was proposed. These schemes exploit the fact that in many applications where (gyro)kinetic models are used, the most stringent CFL condition is associated with the linear part of the model. This observation serves as the basis for the numerical methods we will consider in this work. Starting from the variation of constant formula, the linear part of the model will be solved exactly as part of an exponential integrator, whereas the nonlinear part, which is very often orders of magnitudes less stiff than the linear part, will be treated explicitly in time. In practice, the linear part can then be solved in phase space by using Fourier techniques or semi-Lagrangian schemes and the nonlinear part is approximated by standard finite difference/finite volume/discontinuous Galerkin techniques. 

The numerical results presented in \cite{cep} were generally very favorable. The authors were able to take larger time steps compared to what has been reported for splitting methods in the literature and the computational cost was significantly reduced. In addition, since exponential integrators treat the nonlinear part explicitly, they can be adapted much more easily to different models. Despite these many favorable properties, the largest stable time step size was difficult to predict and varied significantly from method to method. Moreover, as we will see, many exponential integrators can behave rather erratically depending on the specific configuration of the simulation.

Thus, the main goal in this paper is to understand the stability of exponential integrators when applied to purely hyperbolic problems. While there is a large literature and well established theory for exponential integrators applied to parabolic problems (see \cite{ei} and references therein), we will see that for purely hyperbolic problems many surprises are encountered. Based on this analysis we will then propose to use a class of exponential methods, Lawson methods, that do not suffer from the described deficiency. Our analysis explains the efficient and robust behavior of Lawson integrators for this kind of problems. We will also present numerical results for both the Vlasov--Poisson equations and a drift-kinetic model that confirm the expected behavior and shows that using this approach significant performance improvements compared to the exponential integrators used in \cite{cep}, and by extension compared to splitting methods, can be attained.

The paper is organized as follows. First, we offer a brief introduction to exponential methods (section \ref{sec:expint}). Then, in section \ref{ode} a linear stability analysis is performed for both the time and phase space discretization. For the explicit part, we consider both centered differences (such as Arakawa's method) and weighted essentially non-oscillatory schemes (WENO) schemes. In sections \ref{sec:vp} and \ref{sec:dk} we investigate the performance of these methods for the Vlasov--Poisson equations and a four-dimensional drift-kinetic model, respectively.

\section{Exponential integrators and Lawson methods \label{sec:expint}}

Exponential methods are a class of time integration schemes that are applied to differential equations of the form
\begin{equation} \label{eq:expint-eq} \dot{u} = Au + F(u), \end{equation}
where $A$ is a matrix and $F$ is a, in general nonlinear, function of $u$. Usually, both $A$ and $F$ are the result of a spatial semi-discretization of a partial differential equation. Exponential methods are applied to problems where $A$ is stiff or otherwise poses numerical challenges, while $F$ can be treated explicitly.
For the hyperbolic case, a prototypical example is the Vlasov equation \eqref{vlasov}. 
Exponential methods are advantageous if the largest velocity is large compared to the electric field. 
Then the linear part has a much more stringent CFL condition than the nonlinear part of the equation. 
We will consider this example in some detail later in the paper.

In this paper, we will consider two types of exponential methods. The idea of \textit{exponential integrators} is to use the variation of constants formula to rewrite equation \eqref{eq:expint-eq} in the following form
\[ u(t_n+\Delta t) = \exp(\Delta t A) u(t_n) + \int_0^{\Delta t} \exp((\Delta t -s)A) F(u(t_n+s)) \,\mathrm{d}s, \]
where we denote the time step size by $\Delta t >0$ and $t_n = n\Delta t$ with $n\in\mathbb{N}$. This expression is still exact; \ie no approximation has been made. Note, however, that this can not be used as a numerical method as evaluating the integral would require the knowledge of $u(t_n+s)$, which is not available. \textcolor{black}{The idea of an exponential integrator is to approximate the nonlinear part $F(u(t_n+s))$ in terms of the available data. In the simplest case we just evaluate it at the left endpoint. That is, we use $F(u(t_n+s))\approx F(u^n)$. Then we can integrate the term $\exp((\Delta t -s)A)$ exactly and obtain.}
\[ u(t_n+\Delta t) \approx u^{n+1} = \exp(\Delta t A) u^n + \Delta t \varphi_1(\Delta t A) F(u^n), \]
where $\varphi_1(z)=(\mathrm{e}^z-1)/z$ is an entire function. This is the first order exponential Euler method. In a similar way exponential Runge--Kutta methods can be constructed. We refer to the literature, in particular the review article \cite{ei}, for more details.

Another ansatz to remove the stiff linear term from equation \eqref{eq:expint-eq} is to introduce the change of variable
\[ v(t) := \exp(-t A) u(t). \]
Plugging this into equation \eqref{eq:expint-eq} yields
\begin{equation} \label{eq:lawson-eq} \dot{v}(t) = \exp(-t A) F(\exp(t A) v(t)). \end{equation}
Now we apply an explicit Runge--Kutta method to the transformed equation. In the simplest case, applying the explicit Euler scheme yields
\[ v(t_n+\Delta t) \approx v^{n+1} = v^{n} + \Delta t \exp(-t_n A) F(\exp(t_n A)v^{n}). \]
Reversing the change of variables yields
\[ u^{n+1} = \exp(\Delta t A) u^{n} + \Delta t \exp(\Delta t A) F(u^{n}). \]
This is the Lawson--Euler method, also a method of order one. Lawson methods are also commonly referred to as integrating factor methods. We immediately see that any explicit Runge--Kutta method applied to equation \eqref{eq:lawson-eq} uniquely determines a Lawson scheme. We call the chosen Runge--Kutta method the \textit{underlying Runge--Kutta method}. For more details we refer the reader to \cite{lawson1967,canuto1988,trefethen2000,minchev2005}.

The example of the Lawson--Euler method already shows the similarity between Lawson schemes and exponential integrators. In fact, Lawson methods can be considered a subclass of exponential integrators. That is, they are a type of exponential integrators that only involve the exponential, but no other matrix functions. For the purpose of this paper we keep the nomenclature distinct. A Lawson scheme is a numerical method obtained as described above, while an exponential integrator is a numerical scheme that, in addition to the matrix exponential, uses other matrix functions.

The efficiency of exponential methods crucially depends on a good method to evaluate the application of the required matrix functions to a vector. A range of methods has been developed to accomplish this. For example, Krylov methods or interpolation at Leja points can be used for a wide range of problems; see, for example, \cite{higham2008,hochbruck1997,al2011,caliari2014, ckor15}. 
 However, often the most efficient approach is to exploit particular knowledge about the differential equation under consideration. For example, in the hyperbolic case $A$ might be a linear advection operator. In this case the application of $\exp(\Delta t A)$ can be computed by using Fourier techniques or semi-Lagrangian schemes. Much research effort has been dedicated towards improving spectral and semi-Lagrangian schemes for kinetic problems \cite{crouseilles2011,einkemmer2014,einkemmer2016,filbet2003,grandgirard,klimas1994,Morrison2017,rossmanith2011,sonnendrucker1999,cheng1976,sircombe2009valis,crouseilles2015hamiltonian,crouseilles2016asymptotic,einkemmer2019comparison,einkemmer2014convergence,einkemmer2014dG} and obtaining good performance on state of the art HPC systems \cite{rozar2013,einkemmer2015,bigot2013,latu2007gyrokinetic,mehrenberger2013vlasov,einkemmer2016mixed,crouseilles2009parallel,einkemmer2019gpu}.

Before proceeding, let us note that for parabolic problems, \ie~where $A$ is an elliptic operator, a mature theory for exponential integrators is available. We again refer the reader to the review article \cite{ei}. In this setting there are relatively few surprises with respect to stability and even rigorous convergence results are available. In addition, exponential integrators have been considered for problems that include both hyperbolic and parabolic terms (see, for example \cite{martinez2009,tambue2010,einkemmer2016expintmhd,einkemmer2013expintgpu}). An interesting point to make is that in this community Lawson methods have all but lost their appeal. In fact, there are many reasons why exponential integrators are to be preferred. For example, if a Krylov method is used to compute the matrix functions, the $\varphi_1$ function usually converges faster than the exponential. In addition, exponential integrators that retain their full order for non-homogeneous boundary conditions have been constructed~\cite{hochbruck2005}. It has been shown that this property can not be achieved for Lawson methods \cite{hochbruck2017lawson}. However, the situation for purely hyperbolic problems is markedly different. Most of the theoretical results that have been obtained in an abstract framework do not apply and there is relatively little literature available. We will see  that the stability for exponential integrators in the fully hyperbolic setting is full of surprises. Moreover, since for kinetic problems we usually have efficient methods to compute the matrix exponential and complicated boundary conditions are rather rare, Lawson methods are an attractive choice due to their improved stability, as we will see.

\section{Linear analysis \label{ode}}

Determining the stability of a numerical scheme by conducting an analysis of a linear and scalar test equation is very well established in the literature. Usually, the Dahlquist test equation $\dot{u} = \lambda u$ is considered. The justification for this is that a linear ODE can be written as $\dot{u} = A u$. Once the matrix $A$ is diagonalized we essentially obtain the test equation. For linear PDEs we first perform a space discretization. Then the same argument can be applied to the resulting differential equation and stability constraints, such as the famous CFL condition, can be deduced. In the nonlinear case this, of course, only gives an indication for stability. Nevertheless, in many practical problems the theory derived in this fashion agrees very well with what is observed in numerical experiments.

The situation for exponential integrators and Lawson methods is more complicated as we separate two parts of the differential equation. Thus, we will consider the following test equation
\begin{equation}
\label{ode_linear}
\dot{u} = ia u + \lambda u, \qquad a\in\mathbb{R}, \lambda\in\mathbb{C}, \qquad u(0) = u_0\in\mathbb{C}.
\end{equation}
We note that we are here exclusively interested in equations with two hyperbolic parts. The reason why we allow $\lambda$ to lie in the complex plane is that some space discretization schemes introduce numerical diffusion. Thus, the discretization moves the eigenvalues from the imaginary axis to the left half complex plane.

Although this test equation is used frequently in the literature, its use is also frequently criticized. The reason for this criticism is that in the linear case the equation $\dot{u} = Au + Bu$ can only be transformed to the form given in equation \eqref{ode_linear} if $A$ and $B$ are simultaneously diagonalizable. This is a severe restriction which is usually not true in practice. Thus, the test equation, even in the linear case, gives only a necessary condition for stability. While this argument is certainly correct, we emphasize that if a numerical integrator does not work for the test equation~\eqref{ode_linear} there is not much hope that it would work for more complicated problems. Therefore, the test equation is still useful and in fact we will see that many of the deficiencies of exponential integrators observed in practice can be illustrated well using the test equation~\eqref{ode_linear}.


\subsection{Lawson methods}

Applying a Lawson method to the test equation~\eqref{ode_linear} proceeds as follows. First, we introduce the change of variable
\[ v(t)=e^{-iat} u(t). \]
which yields the equation
\[ \dot{v} = e^{-iat} \lambda (e^{iat }v) = \lambda v. \]
Thus, we precisely obtain the  Dahlquist test equation. We now apply an explicit Runge--Kutta method to that equation. It is well known that this results in
\[ v^{n+1} = \phi(z) v^n, \qquad z=\lambda \Delta t, \]
where $\phi$ is the so-called stability function. Reversing the change of variable we obtain
\[ u^{n+1} = e^{ia\Delta t} \phi(z) u^{n}, \qquad z=\lambda \Delta t. \]
The condition for stability is $| e^{ia\Delta t} \phi(z) | = | \phi(z) | \leq 1$. Thus, the linear stability characteristics of a Lawson scheme is identical to that of its underlying Runge--Kutta method. 

This makes the problem rather easy as the stability constraint for explicit Runge--Kutta methods has been extensively studied in the literature. In our present application we are primarily interested in obtaining numerical methods that maximize the part of the imaginary axis that is included in the domain of stability. It is well known that an $s$ stage method can include at most $i [-(s-1),s-1]$. This is, for example, stated as an exercise in \cite[Chapter. IV.2, exercise 3]{sode}. Thus, unfortunately, there is no analog to Runge--Kutta--Chebyshev methods for hyperbolic problems.

\begin{figure}[h]
	\centering
	\includegraphics[width=0.3\textwidth]{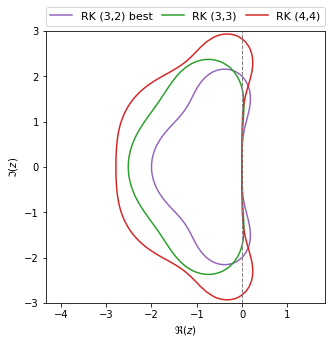}
    \caption{The domain of stability for some classic explicit Runge--Kutta methods is shown. The nomenclature \textit{RK(s,p)} denotes a method with $s$ stages that is of order $p$. The Butcher tableaus of these methods are given in Appendix \ref{butcher}.
    }
	\label{fig:RK_sd2}
\end{figure}

For the sake of completeness we plot in Figure \ref{fig:RK_sd2} the curve given by $|\phi(z)| = 1$ for different Runge-Kutta methods. 
The only non-standard method here is \textit{RK(3,2) best} which is a three stage second order method that has been purposefully constructed to enhance stability on the imaginary axis (see Appendix \ref{butcher} for its Butcher tableau). We also emphasize that the stability domain of the classic four stage fourth order Runge--Kutta method is quite close to the theoretical bound $i [-(s-1),s-1], s=4$.

\subsection{Exponential integrators}

We now apply commonly used exponential integrators to the test equation \eqref{ode_linear}. In this work we will consider the following methods: ExpRK22 (a classic two stage second order method), the method of Cox--Matthews \cite{cox}, the method of Hochbruck--Ostermann  \cite{hochbruck2005}, and the method of Krogstad \cite{krogstad2005}. We refer to \cite{ei} for more details and to Appendix \ref{butcher} for the Butcher tableaus of these methods. 

For the sake of brevity we will only detail the calculation for the ExpRK22 scheme. Applying this method to the test equation we obtain
\begin{eqnarray*}
k_1&=&e^{ia\Delta t}u^n + \Delta t\varphi_1(ia\Delta t)\lambda u^n\nonumber\\
u^{n+1}&=& e^{ia\Delta t}u^n + \Delta t \Big[ (\varphi_1(ia\Delta t)-\varphi_2(ia\Delta t))\lambda u^n + \varphi_2(ia\Delta t)\lambda k_1\Big], 
\end{eqnarray*}
where $\varphi_1(z)=(e^{z}-1)/z$ and $\varphi_2(z)=(e^{z}-1-z)/z^2$ are entire functions. This yields the stability function
\[ \phi(z) = e^{ia\Delta t} + \Big(\varphi_1(ia\Delta t)-\varphi_2(ia\Delta t)+e^{ia\Delta t}\varphi_2(ia\Delta t)\Big)z + \varphi_1(ia\Delta t)\varphi_2(ia\Delta t)z^2, \]
where, as before, we use $z=\lambda \Delta t$.

\begin{figure}[h]
	\centering
	\begin{subfigure}[b]{0.3\textwidth}
		\centering \includegraphics[width=\textwidth]{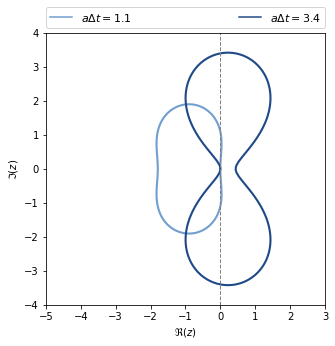}
	\end{subfigure}
	\begin{subfigure}[b]{0.3\textwidth}
		\centering \includegraphics[width=\textwidth]{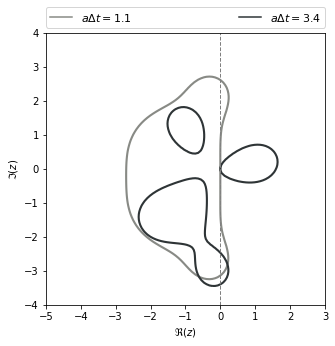}
	\end{subfigure}

	\begin{subfigure}[b]{0.3\textwidth}
		\centering \includegraphics[width=\textwidth]{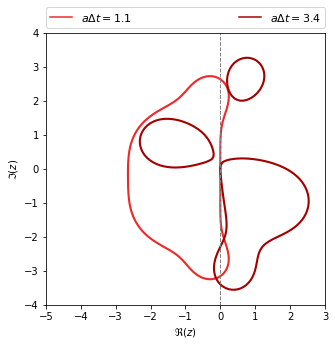}
	\end{subfigure}
		\begin{subfigure}[b]{0.3\textwidth}
		\centering \includegraphics[width=\textwidth]{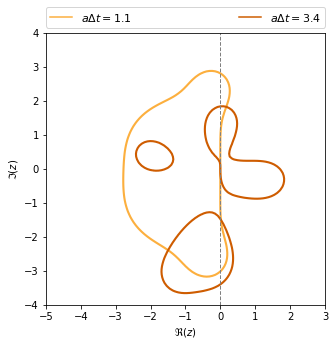}
	\end{subfigure}
    \caption{Stability domain of exponential integrators for two different values of $a\Delta t\in\{1.1, 3.4\}$.  From top left to bottom right: ExpRK22,  Krogstad, Cox--Matthews and Hochbruck--Ostermann.}  
	\label{fig:expRK_sd}
\end{figure}

Our first observation is that, in contrast to Lawson methods, the behavior of this stability function can not be understood by the domain of stability of the underlying $RK(2, 2)$ method, \ie~the explicit method we obtain if we take $a \to 0$. In fact, as we vary $a$ the domain of stability changes drastically. The domain of stability for the four exponential integrators (ExpRK22, Cox--Matthews, Hochbruck--Ostermann, and Krogstad) is plotted in Figure~\ref{fig:expRK_sd} for $a\Delta t=1.1$ and $a\Delta t=3.4$. It is most striking that for large $\vert a\Delta t \vert$ the domain of stability does not contain a symmetric interval of the imaginary axis. It should be evident that this has the potential to causes severe stability issues.

\subsection{Phase space discretization}

We start from the two-dimensional linear transport equation
\begin{equation}
	\label{vp_linear}
    \partial_t f + d\partial_x f + b\partial_v f = 0, \;\; d, b\in\mathbb{R}, \;\; x\in [0, 2 \pi], \;\; v\in [-v_{\max},v_{\max}], 
\end{equation}
where $v_{\max}>0$  refers to the truncated velocity domain. 
The sought-after distribution function is $f(t,x,v)$ and we impose periodic boundary conditions in the $x$-direction. We assume that $d$ and $b$ are constants and thus the corresponding operators commute. This is an idealization of the Vlasov equation we will consider in the next section. In preparation for that example it is most useful to think that $d$ is large and thus would induce a stringent CFL condition if discretized by an explicit scheme.

We now have to discretize this equation both in the $x$ and the $v$ direction. In the spatial direction $x$, we will consider a spectral approximation. Performing a Fourier transformation of equation \eqref{vp_linear} yields
\begin{equation} 
\label{fourier_x_vlasov}
\partial_t \hat{f}_{k} +  i d k\hat{f}_{k} +b \partial_v \hat{f}_{k} = 0,
\end{equation} 
where $\hat{f}_k(t,v)$ denotes the Fourier transform of $f(t,x,v)$ with respect to $x$. The corresponding frequency is denoted by $k$. 

We now perform the discretization in the $v$ direction. The grid points are denoted by $v_j=-v_{\max}+j \Delta v$, with $\Delta v=2v_{\max}/N_v$, where $N_v$ is the number of points. We will consider two options here. Namely, either using a centered difference scheme or an upwind scheme.

\paragraph{Centered scheme in $v$.\\} The classic centered scheme is obtained by approximating the velocity derivative in equation \eqref{fourier_x_vlasov} by
$$
(\partial_v \hat{f}_{k})(v_j) \approx \frac{\hat{f}_{k, j+1}-\hat{f}_{k, j-1}}{2\Delta v},
$$ 
where $\hat{f}_{k,j}$ is an approximation of $\hat{f}_{k}(v_j)$.  Inserting this centered approximation 
in \eqref{fourier_x_vlasov} yields
\begin{equation} \label{eq:half-fourier}
\partial_t \hat{f}_{k,j} +  i d k\hat{f}_{k,j} +b \frac{\hat{f}_{k,j+1} -\hat{f}_{k,j-1} }{2\Delta v} = 0. \end{equation}
The system is already diagonal with respect to the index $k$. We now also diagonalize it with respect to the index $j$. To do that we express the function in terms of its Fourier modes with respect to $v$. That is,
\[ \hat{f}_{k,j} = \sum_m \bar{f}_{k, m}\exp\left(i \frac{2\pi m}{2v_{\max}}  v_j\right), \]
where $\bar{f}_{k,m}$ denotes the (double) Fourier transform of $f$ with frequency in space $k$ and frequency in velocity $m$. Inserting this into equation \eqref{eq:half-fourier} yields
\begin{equation}
\label{discrete_linear_transport}
	\partial_t \bar{f}_{k,m} + i dk\bar{f}_{k,m} +b \frac{i\sin(2\pi m \Delta v/(2v_{\max})) }{\Delta v} \bar{f}_{k,m}= 0.  
\end{equation}
We immediately see that this equation is precisely in the form of equation \eqref{ode_linear} as studied in the previous section. We also observe that $\lambda \in i\mathbb{R}$. That is, the eigenvalues for the centered difference approximation lie exclusively on the imaginary axis. One immediate consequence is that for Lawson methods the CFL condition is given by $b \Delta t < C \Delta v$, where $C$ is chosen such that $i [-C,C]$ lies in the domain of stability of the underlying Runge--Kutta method.

\paragraph{Linearized WENO approximation in $v$.\\} 

A common technique to discretize hyperbolic partial differential equations is to use the so-called weighted essentially non-oscillatory schemes (WENO) schemes. These are nonlinear schemes that limit oscillations in regions where sharp gradients occur, but still yield high order accuracy in smooth regions of the phase space. In the linear case WENO schemes reduce to upwind discretizations. Here, we will consider the LW5 scheme (the linearized version of the WENO5 scheme as considered in \cite{baldauf, lunet, motamed, wang}) that is given by (from now on we assume w.l.o.g.~that $b>0$)
\begin{align*}
(\partial_v \hat{f}_k)(v_j) &\approx \frac{1}{\Delta v}\Big(-\frac{1}{30} \hat{f}_{k,j-3} +\frac{1}{4} \hat{f}_{k,j-2} -\hat{f}_{k,j-1} + \frac{1}{3} \hat{f}_{k,j} +\frac{1}{2} \hat{f}_{k,j+1} - \frac{1}{20} \hat{f}_{k,j+2}\Big). 
\end{align*}
We now perform the same analysis as for the centered scheme (see \cite{baldauf, fov}). This yields
\begin{align}
  {\color{black} \left( \frac{2\pi }{2 v_{\max}}  i m \bar{f}_{k, m} \approx \right)}
   \mu_m \bar{f}_{k, m} &:= \frac{\bar{f}_{k,m}}{\Delta v}\Big(-\frac{1}{30} e^{-\frac{3i m\pi \Delta v}{v_{\max}}} +\frac{1}{4} e^{-\frac{2i m \pi \Delta v}{v_{\max}}} -e^{-\frac{i m \pi \Delta v}{v_{\max}}}  \nonumber\\
&\hspace{1cm}+ \frac{1}{3}  +\frac{1}{2} e^{\frac{i m \pi \Delta v}{v_{\max}}} - \frac{1}{20} e^{\frac{2i m \pi \Delta v}{v_{\max}}}\Big). 
\label{lw5symbol}
\end{align}
We then obtain the equation 
\begin{equation}
	\label{discrete_linear_transport_weno}
	\partial_t \bar{f}_{k,m} + i dk\bar{f}_{k,m} +b \mu_m \bar{f}_{k,m}= 0.
\end{equation}
Once again this is precisely the form of equation \eqref{ode_linear}, where $a=dk\in\mathbb{R}$ and $\lambda=b\mu_m \in\mathbb{C}$. The main difference to the centered difference scheme is that $\lambda$ is not necessarily on the imaginary axis. In fact, the eigenvalues aquire a negative real part which stabilizes the scheme and avoids spurious oscillations, but also adds unphysical dissipation to the numerical method.

{\color{black} 
\begin{remark}
Let us remark that performing a Fourier approximation in $v$ is also possible. Using the same notation as before, the counterpart of \eqref{discrete_linear_transport} and \eqref{discrete_linear_transport_weno} in that case is
$$
\partial_t \bar{f}_{k,m} + i dk\bar{f}_{k,m} +b i \frac{2\pi m}{2 v_{\max}} \bar{f}_{k,m}= 0.
$$
It is worth mentioning that this last equation can be obtained by considering the limit $\Delta v$ goes to zero in  \eqref{discrete_linear_transport} or \eqref{discrete_linear_transport_weno}. The stability condition can be computed as for the centered difference case, since the eigenvalues of the Fourier approximation also lies on the imaginary axis. The CFL condition is given by $b \pi  < C \Delta v$, where $C$ is chosen such that $i [-C,C]$ lies in the domain of stability of the underlying Runge--Kutta method.

\end{remark}
}

\subsection{Computing the CFL condition}

Equipped with the knowledge of the domain of stability for the time discretization and the eigenvalues of the space discretization, we are now in a position to determine the CFL condition for the linear transport equation \eqref{vp_linear}. This task will be rather easy to accomplish for the Lawson schemes, where the stability does not depend on the advection speed for the transport in the $x$ direction. However, for exponential integrators even this linear analysis is rather complicated, as we will see.

\subsubsection{Centered scheme in $v$.}
In the case of centered approximation of the velocity derivative, the Fourier multiplier is a 
pure imaginary complex number (see equation \eqref{discrete_linear_transport}). We thus look for $y_{\max}\in \mathbb{R}_+$ such that the interval $i(-y_{\max},y_{\max}) \subset {\cal D}$, where ${\cal D}$ is the domain of stability for the chosen time integrator. 

\paragraph{Lawson integrators.\\} 
We simply look for the largest value $y_{\max}$ such that $i(-y_{\max}, y_{\max}) \subset \mathcal{D}$. The corresponding values for a number of schemes are listed in Table \ref{tab:ymax_Lawson}. These values have to be understood in the 
following way: they induce the CFL condition $b \Delta t\leq y_{\max}\Delta v$ for the discretized equation \eqref{discrete_linear_transport}, where $\Delta t$ denotes the time step size and $\Delta v$ is the velocity mesh size.

\begin{table}[h]
	\centering
	\begin{tabular}{|c|c|c|c|}
		\hline
		Methods & Lawson($RK(3,2) \; best$) & Lawson($RK(3,3)$) & Lawson($RK(4,4)$) \\
		\hline
		$y_{\max}$ & $2$ & $\sqrt{3}$ & $2\sqrt{2}$\\
		\hline  
	\end{tabular}
	\caption{CFL number for some Lawson schemes applied to \eqref{discrete_linear_transport}. }
	\label{tab:ymax_Lawson}
\end{table}

\paragraph{Exponential integrators. \\}

For the exponential integrators the domain of stability is very sensitive to the value of $(a\Delta t)$. To get an idea of what we can expect, we consider the quantity $y_{\max} = \min_{(a\Delta t)\in\mathbb{R}} \; y^{exp}_{\max}(a \Delta t)$. As before, $y^{exp}_{\max}(a\Delta t)$ is the largest value such that $i(-y^{exp}_{\max}(a\Delta t), y^{exp}_{\max}(a\Delta t))\subset {\cal D}$, where ${\cal D}$ is the domain of stability for the chosen exponential time integrator for a given $(a\Delta t)$. 
Even for relatively simple numerical methods it is not possible to compute this quantity analytically. Thus, we resort to numerical approximations. Unfortunately, it turns out that for most exponential integrators this value is zero. This can be appreciated by considering Figure \ref{fig:expRK_sd} once more. Clearly, there are values of $(a \Delta t)$ such that no relevant part of the imaginary axis (or only half the imaginary axis) is part of the domain of stability.
Thus, most exponential integrators are unstable in the von Neumann sense. However, this is not what we observe in practice. In fact, already the results presented in \cite{cep} indicate that we can successfully run numerical simulations using, for example, the Cox--Matthews scheme. There are two major points to consider here
\begin{itemize}
    \item The $y_{\max}$ obtained is a worst case estimate. In fact, we know that for $\Delta t \to 0$ we regain the stability of the underlying Runge--Kutta method. Thus, for small $(a \Delta t)$ the methods is expected to work well.
    \item As is usually done we have mandated that $\vert \phi(z) \vert \leq 1$. However, strictly speaking this is not necessary for practical simulation. If we assume that $\vert \phi(z) \vert \leq 1+\varepsilon$ and we take $n$ steps the amplification of the error is given by $(1+\varepsilon)^n$. In the limit $n \to +\infty$ this quantity diverges. However, since we usually do not take infinitely small time steps we still can hope to obtain a relatively accurate approximation, especially if $\varepsilon$ is small. \textcolor{black}{In particular, if $\varepsilon = C \Delta t = C t_{final}/n$ (where $t_{final}$ denotes the final time and $n$ the number of iterations) we have $(1+\varepsilon)^n = (1+\tfrac{C t_{final}}{n} )^n \leq \exp(C t_{final})$ and thus the scheme is stable, while the error constant is increased by $\exp(C t_{final})$.}
\end{itemize}

To investigate this further, we propose to relax the stability condition by introducing a threshold $\varepsilon>0$ 
in the definition of the stability domain
\begin{equation}
\label{d_eps}
	\mathcal{D}_\varepsilon = \{ z\in\mathbb{C} : |\phi(z)| \leq 1+\varepsilon \}. 
\end{equation}

\begin{figure}[h]
  \centering
  \begin{subfigure}[b]{0.33\textwidth}
        \centering \includegraphics[width=\textwidth]{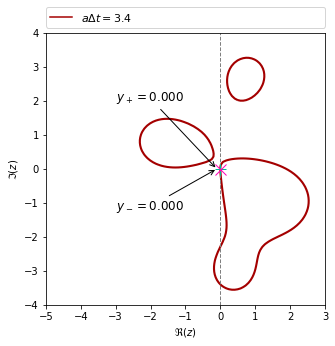}
  \end{subfigure}
  \begin{subfigure}[b]{0.33\textwidth}
        \centering \includegraphics[width=\textwidth]{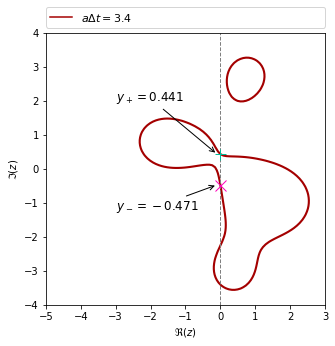}
  \end{subfigure}
  \caption{Example of variation of $y_+$ and $y_-$ when we relax the stability condition for the Cox--Matthews scheme. We represent $\mathcal{D}_{0}$ on the left, and $\mathcal{D}_{10^{-2}}$ on the right with the values $y_+$ and $y_-$  such that 
  $i(y_-, y_+)\subset \mathcal{D}_{\varepsilon}$.}
  \label{ymax_example}
\end{figure}

In Figure~\ref{ymax_example}, we plot the domain of stability for the Cox--Matthews method and $a\Delta t = 3.4$ for $\varepsilon=0$ and $\varepsilon=10^{-2}$. One can observe that in the latter case a non-zero $y_{\max}^{exp}(3.4)$ is obtained. We also call attention to the fact that the part of the imaginary axis included in this relaxed stability domain is not symmetric. 

\begin{figure}[h]
	\centering
	\includegraphics[scale=0.3]{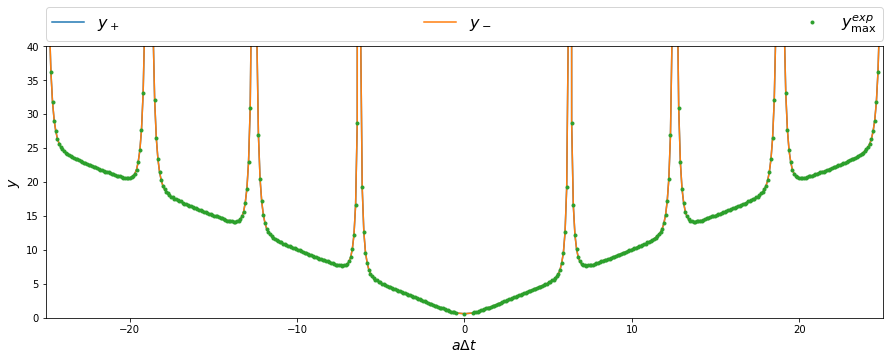}
    \caption{$y^{exp}_{\max}$, $|y_+|$ and $|y_-|$ as a function of $a \Delta t$ for the ExpRK22 method with $\varepsilon=10^{-2}$. } 
	\label{ymax_expRK22}
\end{figure}

In Figure \ref{ymax_expRK22}, we plot the dependence of $y^{exp}_{\max}$ as a function of $(a\Delta t)$ for $\varepsilon=10^{-2}$ and the ExpRK22 method. Let us recall that for $\varepsilon=0$, the method gives $y_{\max}=0$. One can observe that the  domain of stability 
${\cal D}_\varepsilon$  of this method is still symmetric with respect to the real axis. In addition, the method becomes more stable as $|a\Delta t|$ increases. Thus, the behavior of the method is completely different from the configuration with $\varepsilon=0$.

\begin{figure}[h]
	\centering
	\includegraphics[scale=0.3]{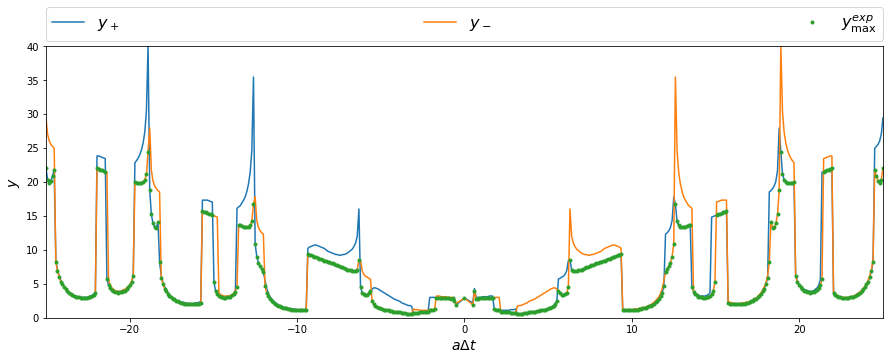}
	\caption{$y^{exp}_{\max}$, $|y_+|$ and $|y_-|$  as a function of $a\Delta t$ for the Hochbruck--Ostermann method with $\varepsilon=10^{-2}$.} 
	\label{ymax_HO}
\end{figure}

In Figure \ref{ymax_HO}, we plot $y^{exp}_{\max}$ as a function of $(a\Delta t)$ for the Hochbruck--Ostermann method (once again for 
$\varepsilon=10^{-2}$). This schemes also gives $y_{\max}=0$ for $\varepsilon=0$. In this case the domain of stability is not symmetric 
and the stability depends quite erratically on the value of $(a \Delta t)$.
In Table \ref{tab:ymax_expo} we have summarized the values of $y_{\max}$ for the four exponential integrators considered in this paper.

\begin{table}
	\centering
	\begin{tabular}{|c|c|c|c|c|}
		\hline
		Methods                                & ExpRK22 & Krogstad & Cox--Matthews & Hochbruck--Ostermann      \\
		\hline
		$y^{exp}_{\max} (\varepsilon=10^{-3})$ & $0.300$ & $0.100$  & $0.150$      & $0.250$ \\
		\hline
		$y^{exp}_{\max} (\varepsilon=10^{-2})$ & $0.551$ & $0.200$  & $0.450$      & $0.501$ \\
		\hline  
		$y^{exp}_{\max} (\varepsilon=10^{-1})$ & $1.001$ & $0.601$  & $1.351$      & $1.702$ \\
		\hline  
	\end{tabular}
	\caption{CFL number, assuming the relaxed stability constraint, for some exponential integrators applied to \eqref{discrete_linear_transport}.}
	\label{tab:ymax_expo}
\end{table}

{\color{red}
Finally, we give an illustration of the above comments by running \eqref{vp_linear} ($d=b=1$ and $v_{\max}=3$) 
with a discontinuous initial data to study the impact of high-frequency on the stability of exponential schemes coupled 
with a centered scheme in $v$. 
The initial condition is 
$$
f_0(x, v) = 1 \mbox{ if } \sqrt{(x-\pi)^2+v^2} \leq 1, \mbox{ and } 0  \mbox{ elsewhere}. 
$$
From the stability analysis we know that for each $k$ it holds that $\bar{f}^{n+1}_{k,m} = \phi(z)\bar{f}^{n}_{k,m}$,
where $\phi$ is the stability function. 
Our study enables us to estimate the amplification factor $|\phi(z)|$ by $(1+\varepsilon)$ (uniformly with respect to $k$)  
so that $|\bar{f}^{n+1}_{k,m}|^2 \leq  (1+\varepsilon)^2 |\bar{f}^{n}_{k,m}|^2$. Hence, for a given $\varepsilon$, 
we consider a time step according to Table \ref{tab:ymax_expo} and run 
two exponential schemes, namely Hochbruck-Ostermann and Cox-Matthews,   
for $100$ timesteps. After each time step
we compute $\| f^n\|^2_{\ell^2} / \| f^0\|^2_{\ell^2}$. 
The results can be found in Figure \ref{instab}.   
\begin{figure}
\centering
\begin{tabular}{cc}
\includegraphics[scale=0.5]{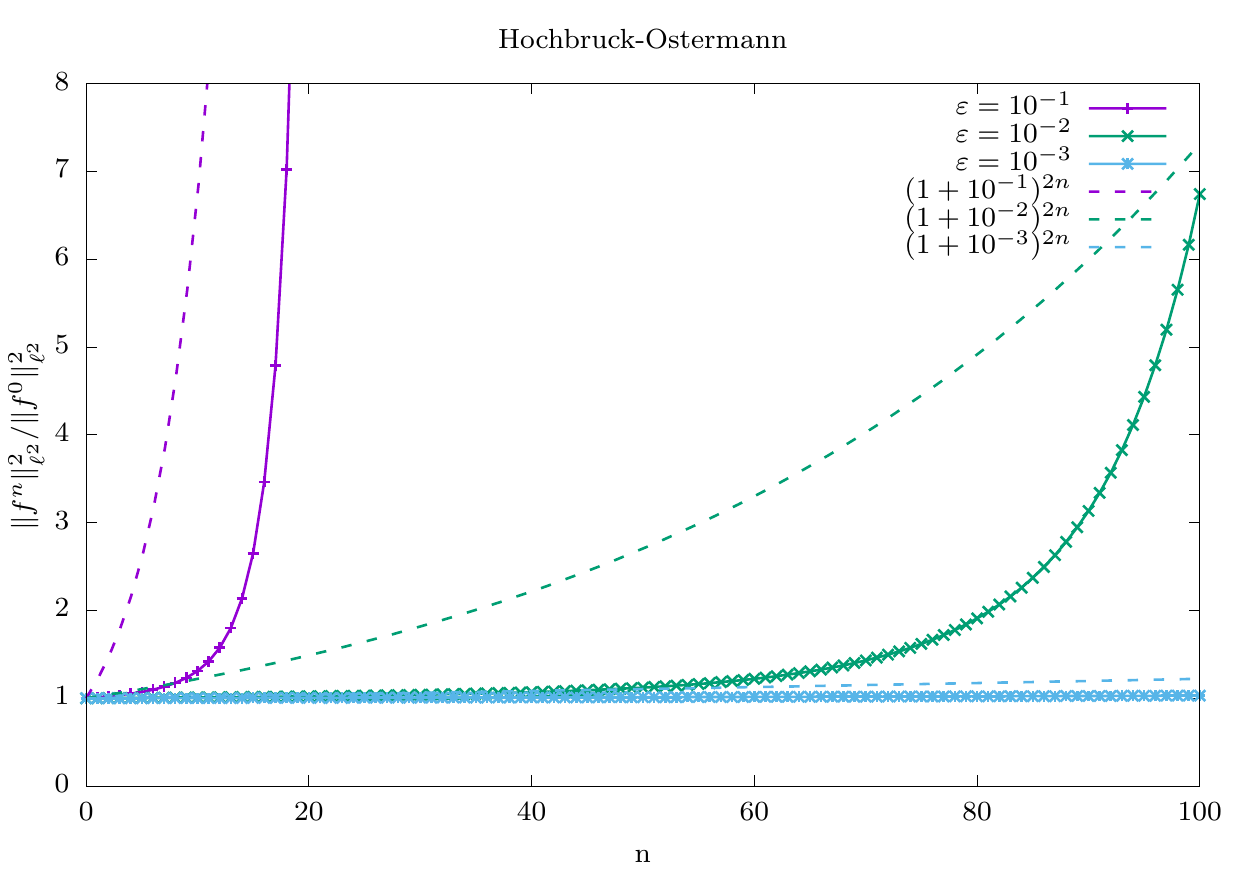} & 
\includegraphics[scale=0.5]{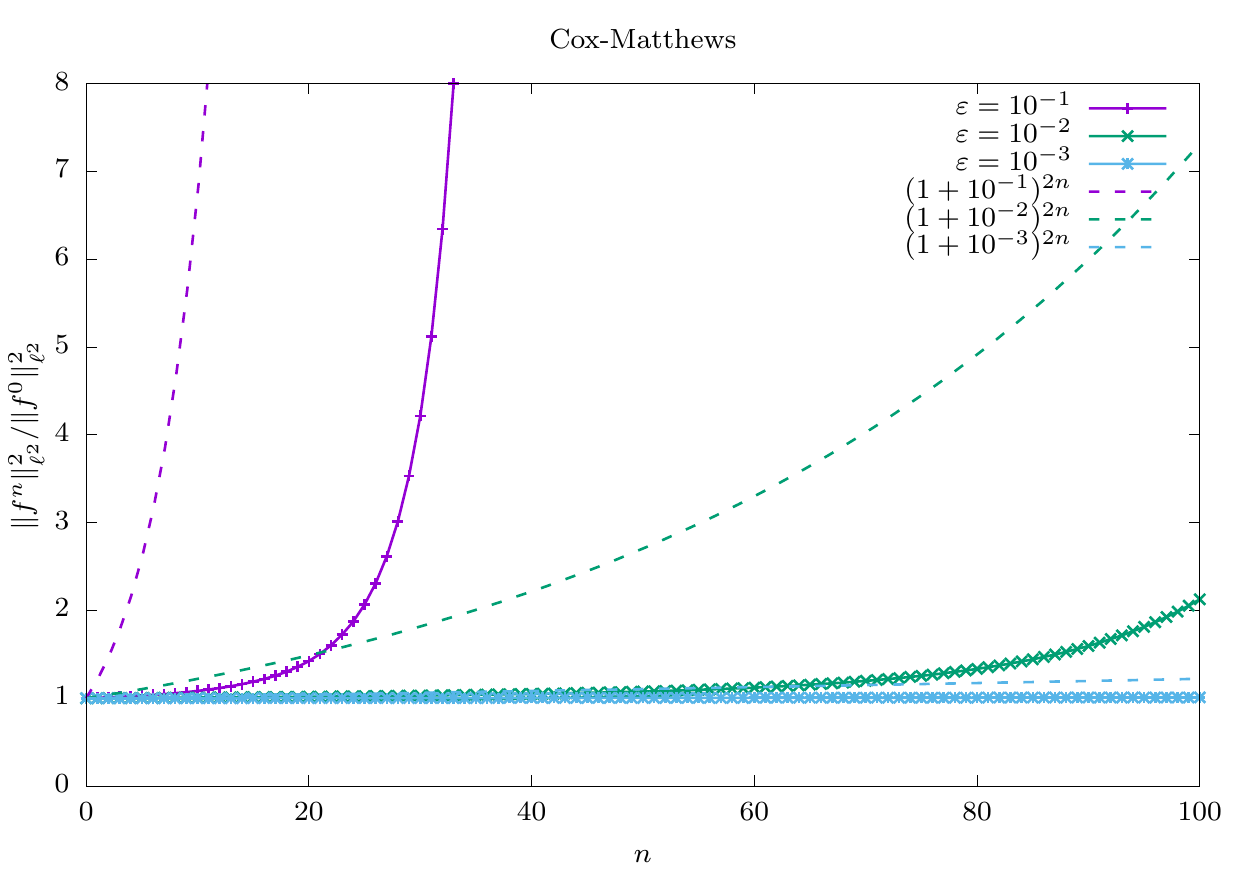} 
\end{tabular}
\caption{Evolution of $\| f^n\|^2_{\ell^2} / \| f^0\|^2_{\ell^2}$ and of $(1+\varepsilon)^{2n}$ as a function of $n$ for different values of 
$\varepsilon$. Left: Hochbruck-Ostermann scheme. Right: Cox-Matthews scheme. }
\label{instab}
\end{figure}
First, we have a numerical confirmation that $(1+\varepsilon)$ is an estimate of the amplification factor 
(the ratio $\| f^n\|^2_{\ell^2} / \| f^0\|^2_{\ell^2}$  always lies under the curve $n\mapsto (1+\varepsilon)^{2n}$).  
Second, since the number of time steps is fixed, for $\varepsilon=0.1$ the simulation becomes unstable for $n\geq 20$. 
As soon as $\varepsilon$ is small enough, the simulation is stable. 
}

\subsubsection{Linearized WENO5 (LW5) scheme. }

In the case of a WENO5 approximation of the velocity derivative, we can not 
easily find a Fourier multiplier because of its nonlinearity. 
Recent studies about stability of WENO5 \cite{wang, motamed, lunet} considered 
the linearized version of WENO schemes by freezing the nonlinear weights, 
so that WENO5 reduces to a high order (linear) upwind scheme called LW5. For this LW5 scheme we can
compute the eigenvalues, see equation \eqref{lw5symbol}, 
and different time stepping schemes can be studied to determine the stability limit. 
We consider only Lawson methods here since we found 
that the exponential schemes we considered (\ie{}~ExpRK22, Krogstad, Cox--Matthews, Hochbruck--Ostermann) 
lead to unstable results when they are combined with LW5 (even in the weak sense considered in the previous section).

The goal  is then to determine the largest non-negative real number $\sigma>0$ such that the 
eigenvalues of the upwind scheme LW5 scaled by $\sigma$ are contained in the domain of stability for the 
time integrator. Since the eigenvalues of LW5 are not as simple as in the case of the centered scheme, we determine $\sigma$ numerically. The main idea of the algorithm to obtain an estimation of $\sigma$ is:
\begin{enumerate}
    \item The argument $\varphi$ of the eigenvalues $\mu_m$  (normalized by $\Delta v$) given by \eqref{lw5symbol} are discretized using a fine angular grid $\{ \varphi_k\} \subset [-\pi/2, \pi/2]$, since the real part of $\mu_{m}$ is negative due to its diffusive character.
  \item A discretized version of the boundary of the stability domain of the underlying Runge--Kutta method is computed.
  \item For each discretized eigenvalue, we look for the closest boundary point of the Runge--Kutta stability domain. 
  This enables us to compute the associated stretching factor $\sigma(\varphi_k)$.
  \item Taking the minimum over all the discretized eigenvalues yields $\sigma:=\min_{k} \sigma(\varphi_k)$. 
\end{enumerate}

\begin{figure}[h]
  \centering
  \begin{subfigure}[b]{0.4\textwidth}
        \centering \includegraphics[width=\textwidth]{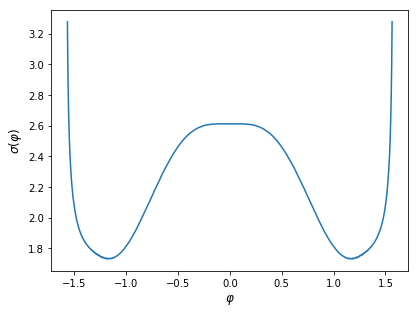}
  \end{subfigure}
  \begin{subfigure}[b]{0.3\textwidth}
        \centering \includegraphics[width=\textwidth]{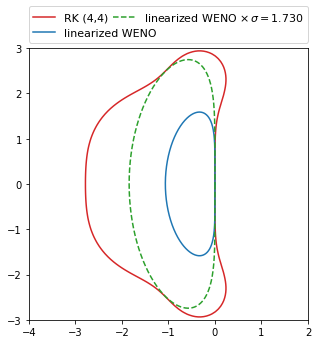}
  \end{subfigure}
  \caption{Left: $\sigma$ as a function of the angle $\varphi$. Right: stability domain of Lawson($RK(4,4)$) (red), 
  eigenvalues for LW5 normalized by $\Delta v$ (blue) and eigenvalues for LW5 normalized by $\Delta v$ 
  stretched with factor $\sigma=1.73$ (dashed green). } 
  \label{cfl_rk44_lw5}
\end{figure}

In Figure \ref{cfl_rk44_lw5} (left), we plot the dependence of $\sigma$ with respect to the angle $\varphi\in [-\pi/2, \pi/2]$ 
for Lawson($RK(4, 4)$) coupled with LW5. We also plot (Figure \ref{cfl_rk44_lw5} (right)) the stability domain of Lawson($RK(4, 4)$), 
the eigenvalues for LW5 (normalized by $\Delta v$) and the eigenvalues for LW5 scaled by $\sigma$. 
The CFL number for some Lawson schemes is shown in Table \ref{tab:ymax_weno_Lawson}. It is interesting to note that the time step size is  reduced compared to the centered scheme.
\begin{table}[h]
	\centering
	\begin{tabular}{|c|c|c|c|}
		\hline
		Methods & Lawson($RK(3,2) \; best$) & Lawson($RK(3,3)$) & Lawson($RK(4,4)$) \\
		\hline
		$\sigma $ & $1.344$ & $1.433$   & $1.73$   \\
		\hline  
	\end{tabular}
	\caption{CFL number for some Lawson schemes applied to \eqref{discrete_linear_transport_weno}. }
	\label{tab:ymax_weno_Lawson}
\end{table}

\section{Numerical simulation: Vlasov-Poisson equations \label{sec:vp}}

In this section we apply Lawson methods and exponential integrators to the Vlasov-Poisson system. We will see that the linear theory developed in the last section gives a good indication of the stability even for this nonlinear problem. We consider a distribution function $f(t, x, v)$ depending on time $t\geq 0$, space $x$, with periodic boundary conditions, and velocity $v$, which satisfies the Vlasov equation 
\begin{equation}
	\label{vlasov}
    \partial_t f(t, x, v) + v\partial_x f(t, x, v) + E(f)(t, x)\partial_v f(t, x, v) = 0   
\end{equation}
coupled to a Poisson problem for the electric field $E(f)(t, x)$ 
\begin{equation}
    \partial_x  E(f)(t, x)= \int_{\mathbb{R}} f(t, x, v) \,dv -1. 
\end{equation}
We employ a Fourier approximation in space. In velocity we either use a centered discretization
$$
\partial_t \hat{f}_{k,j} + v_j i k \hat{f}_{k,j} + \Widehat{\Big(E_{\cdot} \frac{f_{\cdot,j+1}  -f_{\cdot,j-1}}{2\Delta v}\Big)}_k = 0
$$  
or the WENO5 discretization
\begin{equation}
    \displaystyle\partial_t \hat{f}_{k,j} + v_j ik \hat{f}_{k,j} + \Widehat{\Big(E^+_{\cdot} \frac{f^+_{\cdot,j+1/2} - f^+_{\cdot,j-1/2}}{\Delta v} \Big)}_k 
+ \Widehat{\Big(E^-_{\cdot} \frac{f^-_{\cdot,j+1/2} - f^-_{\cdot,j-1/2}}{\Delta v} \Big)}_k = 0, 
\label{vp_weno}
\end{equation} 
where $E^+=\max(E, 0)$, $E^-=\min(E, 0)$ and $f^\pm_{j+1/2}$ denote the numerical fluxes (see Appendix \ref{app_weno} for more details). Both of these phase space discretizations can be easily cast in the following form
\[ \partial_t \hat{f}_{k,j} = - v_j i k \hat{f}_{k,j} + F(f)_{k,j}, \]
for an appropriately defined $F$. We can now apply an exponential integrator or a Lawson scheme. To illustrate this let us consider the exponential Euler method. This gives
\[ \hat{f}^{n+1}_{k,j} = \exp(-\Delta t v_j i k) \hat{f}^n_{k,j} + \Delta t \varphi_1(-\Delta t v_j i k) F(f^n)_{k,j}. \]
Since in Fourier space the exponential and $\varphi_1$ functions have only scalar arguments, their computation is easy and efficient (\ie~no matrix functions have to be computed). Due to the nonlinearity, it is favorable to compute $E \partial_v f$ in real space. This is done efficiently by using the fast Fourier transform. Generalizing this scheme to multiple dimensions in both space and velocity is straightforward.

To apply our theory from the linear analysis to the nonlinear Vlasov-Poisson case, we need a way to compute the CFL condition. Note that the coefficient of the $v$ advection depends on $E$ and thus implicitly on time. We choose the time step for the centered scheme as follows
\begin{equation}
\label{cfl_vlasov_lc}
\Delta t_n = \frac{y_{\max} \Delta v}{\| E^n \|_{L^\infty}}, 
\end{equation}
whereas for the WENO5 scheme we use the CFL condition computed from its linearized version (LW5)
\begin{equation}
\label{cfl_vlasov_lw}
\Delta t_n = \frac{\sigma \Delta v}{\| E^n \|_{L^\infty}}.
\end{equation}
The value $\| E^n \|_{L^\infty}$ is just the maximal value of the electric field at time $t_n$. The values for  $y_{\max}$ and $\sigma$ are given in Tables \ref{tab:ymax_Lawson}, \ref{tab:ymax_expo}, and \ref{tab:ymax_weno_Lawson} according to the chosen time integrator.

\subsection{Landau damping test.} 

We present numerical results for the standard Landau damping test case. The initial condition is given by
$$
f_0(x, v) = \frac{1}{\sqrt{2\pi}} e^{-\frac{v^2}{2}} (1+0.001 \cos(0.5 x)), \;\; x\in [0, 4\pi], v\in \mathbb{R}. 
$$
The numerical parameters are chosen as follows: the number of points in space is $N_x=81$ whereas the velocity domain 
is truncated to $[-v_{\max}, v_{\max}]$ with $v_{\max}=8$ and is discretized with $N_v=128$ grid points.

Let us remark that for the Landau damping test, the conditions \eqref{cfl_vlasov_lc} and 
\eqref{cfl_vlasov_lw} allow us to take very large time steps, since $\| E^n \|_{L^\infty} \leq \| E^0 \|_{L^\infty} = 2\cdot 10^{-3}$. Then, 
we get $\Delta t = C \Delta v \; 0.5\cdot 10^3 =  62.5 C$, where $C$ can be either $y_{\max}$ or $\sigma$ depending on the chosen time integrator. This means that in practice we can choose the time step $\Delta t$ independently from the mesh. This is clearly a desirable desirable feature of the time integrator.

\begin{figure}[h]
	\centering
\includegraphics[width=\textwidth]{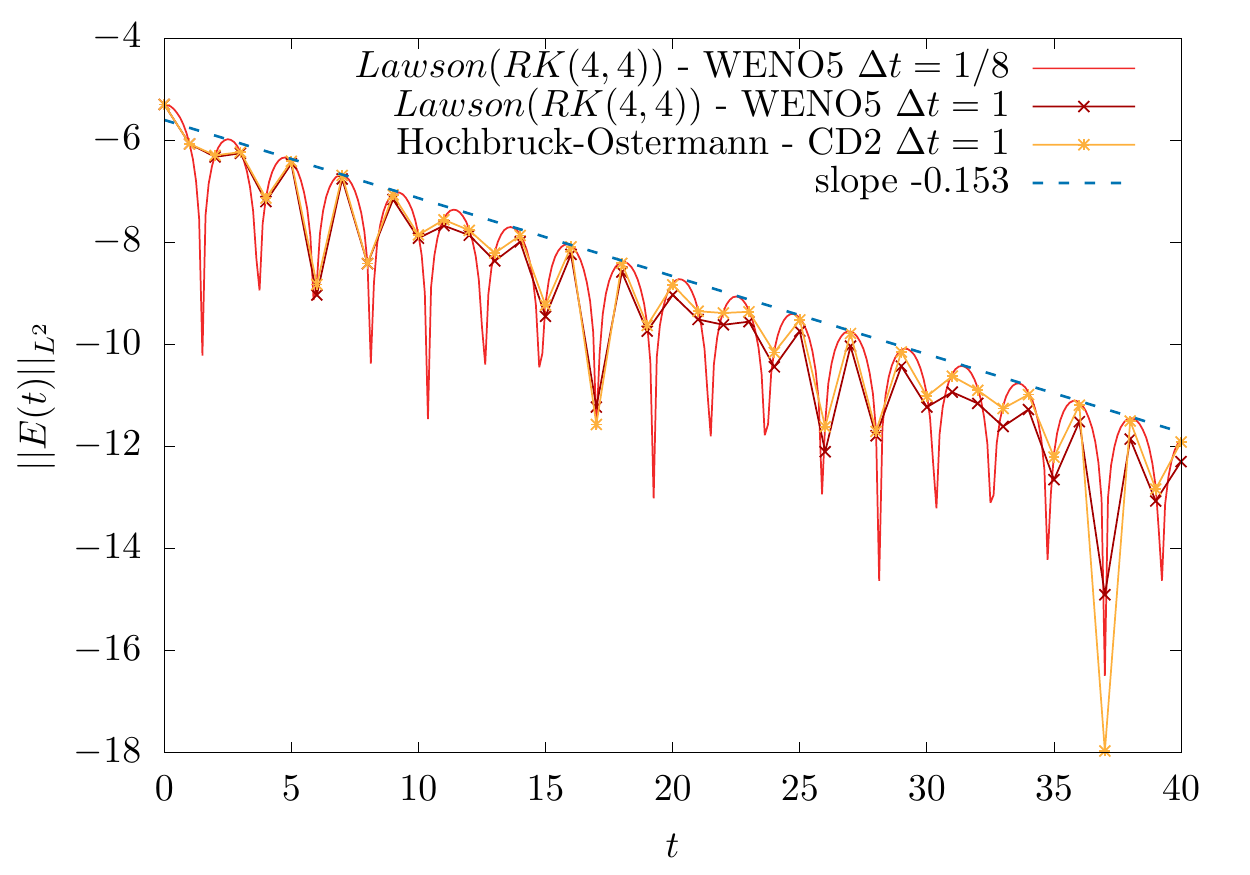}
	\caption{Landau damping test: time history of $\|E(t)\|_{L^2}$ (semi-log scale) obtained with Lawson($RK(4, 4)$) and WENO5 
	(with $\Delta t=1/8$ and $\Delta t=1$) {\color{black} and with Hochbruck-Ostermann and CD2 (with $\Delta t=1$)}.}
	\label{ld}
\end{figure}

In Figure \ref{ld}, the time history of the electric energy $\|E^n\|_{L^2}$ (in semi-log scale) 
{\color{black} using Lawson($RK(4, 4)$)-WENO5 (with two different time steps $\Delta t=1/8$ and $\Delta t=1$), 
and using Hochbruck--Ostermann-CD2 (with $\Delta t=1$).   
One can observe that the expected damping rate  
($\gamma=-0.153$) is recovered for the three schemes. Although, the accuracy deteriorates for $\Delta t=1$ Lawson($RK(4, 4)$)-WENO5, 
the Hochbruck--Ostermann-CD2 scheme gives very good results even with $\Delta t=1$. We note that all the numerical schemes are clearly stable. }

\subsection{Bump on tail test.}

Next, we consider the bump on tail test for which the initial condition is 
$$
f_0(x, v) = \left[\frac{0.9}{\sqrt{2\pi}} e^{-\frac{v^2}{2}} + \frac{0.2}{\sqrt{2\pi}} e^{-2(v-4.5)^2} \right](1+0.04 \cos(0.3 x)), \;\; x\in [0, 20\pi], v\in \mathbb{R}. 
$$
The numerical parameters are chosen as follows: the number of points in space is $N_x=135$ 
whereas the velocity domain $[-v_{\max}, v_{\max}]$ (with $v_{\max}=8$) is discretized with $N_v=256$ grid points. 
Concerning the time step, 
as in the Landau damping example, the conditions \eqref{cfl_vlasov_lc} and \eqref{cfl_vlasov_lw} turn out to be very light for Lawson schemes. 
Indeed, we found $\max_n \| E^n \|_{L^\infty} \approx 0.6$ 
so that, with the considered velocity grid, the time step has to be smaller than $0.14$ 
in the worst case (Lawson($RK(3, 2) \; best)$ combined with WENO5). 
To capture correctly the phenomena involved in the bump on tail test, we take the following time step size
\begin{equation}
\label{dtbot}
\Delta t_n = \min \Big( 0.1,  \frac{C \Delta v}{\|E^n\|_{L^\infty}} \Big), 
\end{equation}
with $C=y_{\max}$ or $\sigma$ depending on the chosen scheme. Thus, also in this configuration we are mostly limited by the accuracy 
and not by the stability constraint.

\begin{figure}
\begin{tabular}{ccc}
    \includegraphics[width=\textwidth]{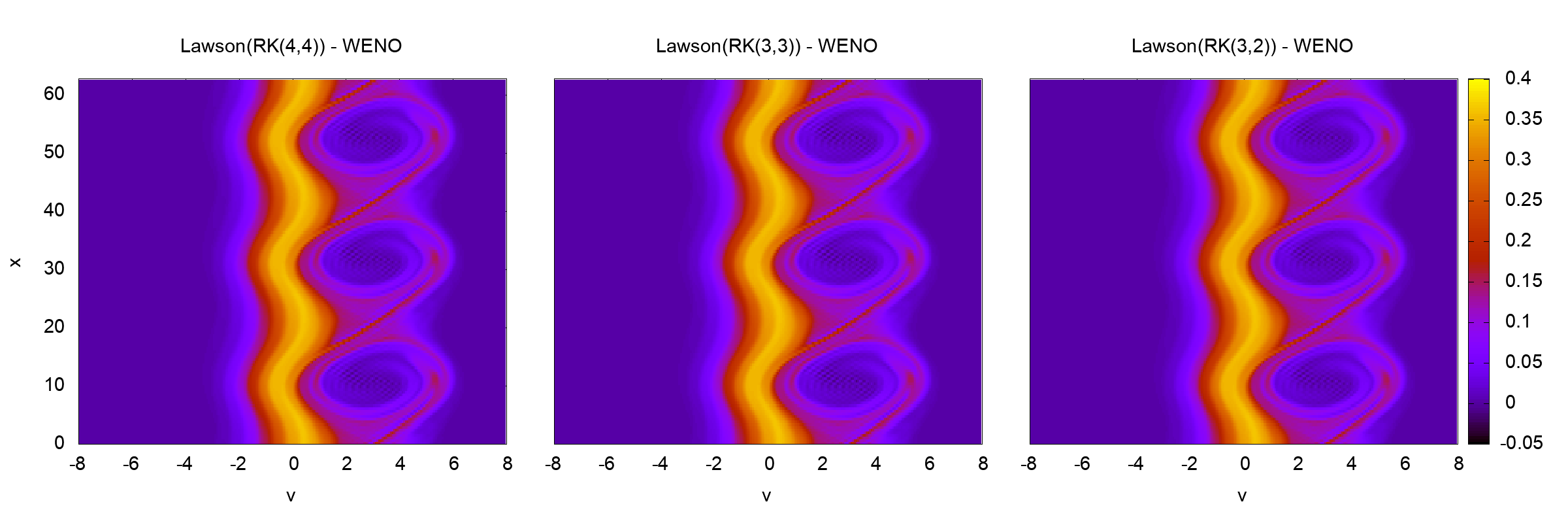} \\
    \includegraphics[width=\textwidth]{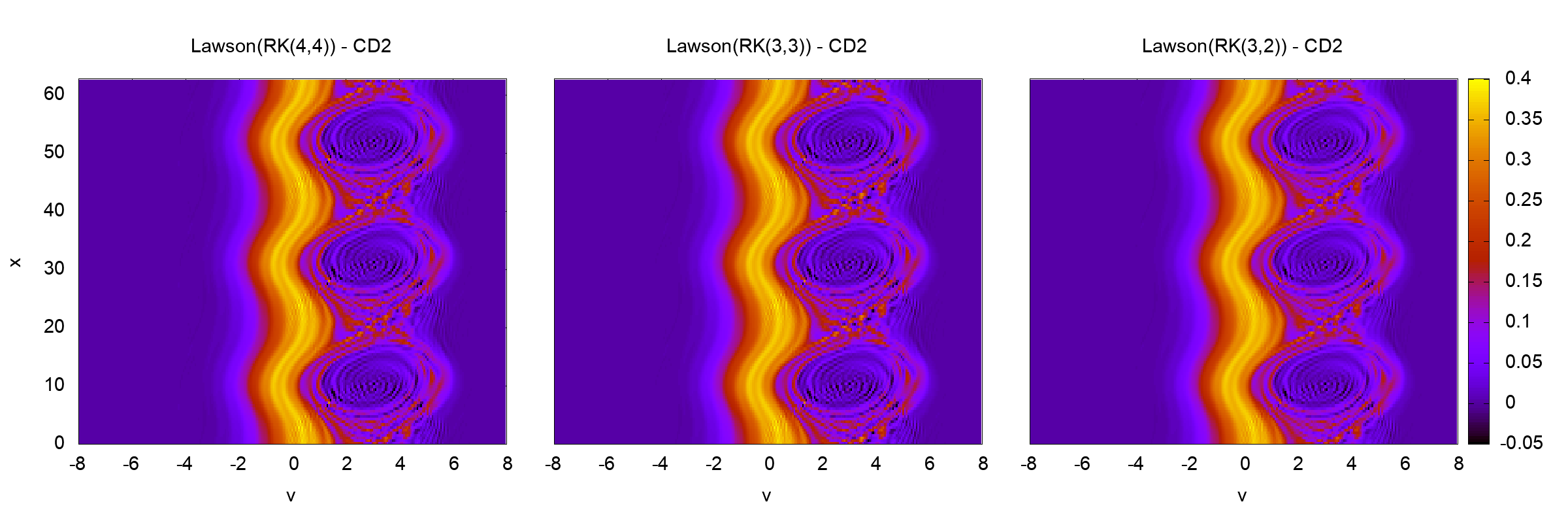}\\
    \includegraphics[width=\textwidth]{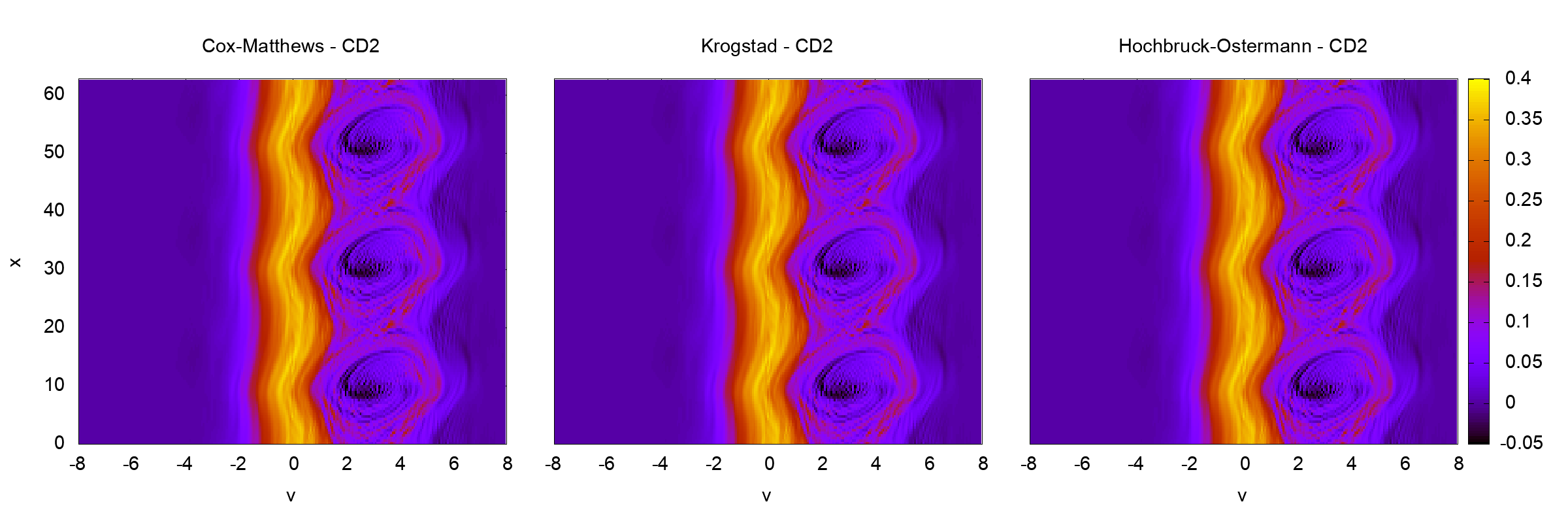}
  \end{tabular}
    \caption{{\color{black} Distribution function at time $t=40$ as a function of $x$ and $v$ for: $(i)$ Lawson schemes ($RK(4, 4)$, $RK(3, 3)$, $RK(3, 2)$) + WENO5 (first row) ; $(ii)$ Lawson schemes ($RK(4, 4)$, $RK(3, 3)$, $RK(3, 2)$) + centered difference scheme (second row) ; $(iii)$ exponential schemes (Cox-Matthews, Krogstad, Hochbruck--Ostermann) + centered difference scheme (third row).}}  
\label{space}      
\end{figure}

{\color{black} In Figure \ref{space}, the full distribution function $f$ is plotted at time $t=40$ ($\Delta t=0.05$) for different schemes (exponential or Lawson in time 
and WENO or centered differences in velocity). One can observe spurious oscillations when the centered differences scheme case is used (second and third rows) whereas the slope 
limiters of WENO5 (first line) are able to control this phenomena so that extremas are well preserved. This is consistent with what has been observed in the literature.}


\begin{figure}
	\centering
  \includegraphics[width=0.9\textwidth]{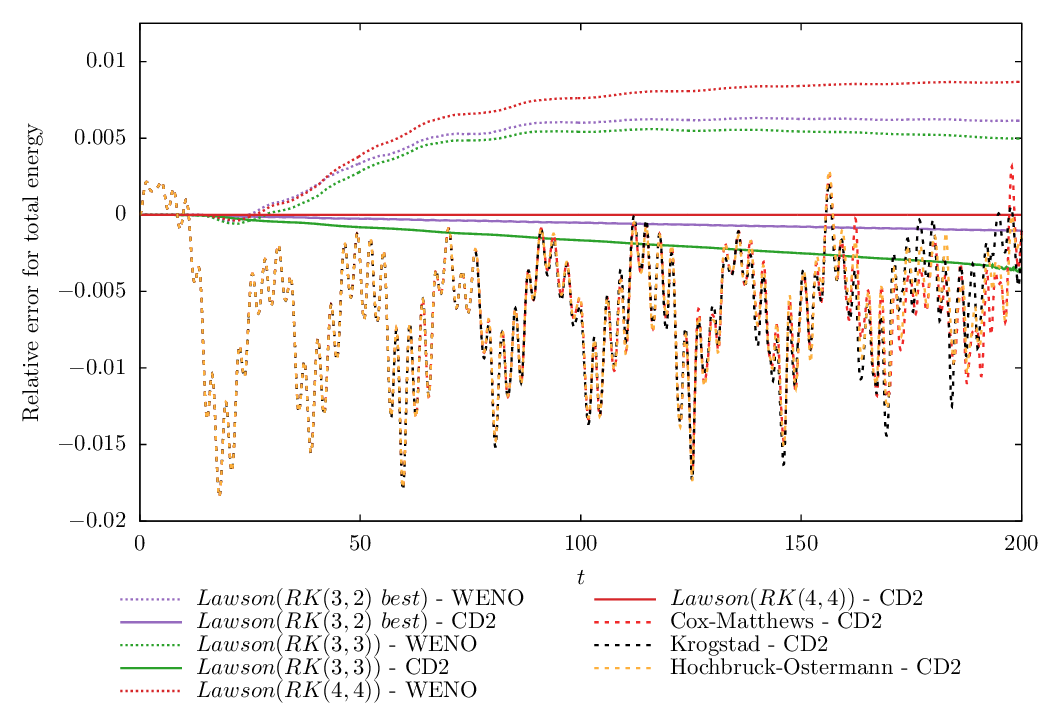}
	\caption{Time evolution of the relative error of the total energy for the different methods. }
	\label{total_energy}
\end{figure}

In Figure \ref{total_energy}, we plot the time evolution of $({\cal H}^n -{\cal H}(0))/{\cal H}(0)$, where ${\cal H}^n\approx {\cal H}(n\Delta t)$ and ${\cal H}(t)$ is the total energy defined by 
$$
{\cal H}(t) = \frac{1}{2}\int\int |v|^2 f(t, x, v) \,dxdv + \frac{1}{2}\int |E|^2(t, x) \,dx. 
$$
This quantity is known to be preserved with time {\color{black} at the continuous level. It is thus a useful metric to evaluate and compare the different numerical 
methods. At this stage, all the used numerical methods are stable and we now look at their accuracy with respect to conservation of energy. }
We observe that Lawson/centered schemes (referred as 'CD2' in the 
legend) preserve this quantity well. It is well known (see for example \cite{cf}) that centered schemes are better at preserving the total energy compared to upwind schemes. The reason is that upwind schemes introduce numerical diffusion. The exponential integrators that are considered show all very similar behavior with respect to energy conservation. They seem to include less drift than the Lawson methods, but for the time scales considered here their error is larger.

\begin{figure}
\centering
\begin{tabular}{cccc}
  \includegraphics[scale=0.3]{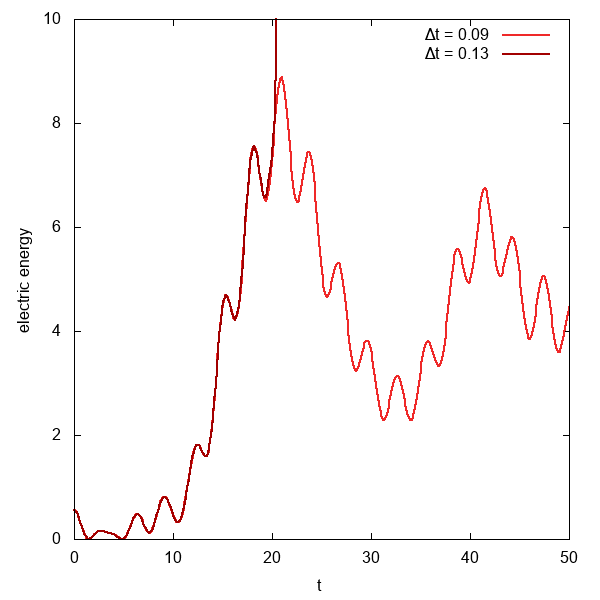} &\hspace{-0.2cm}\includegraphics[scale=0.3]{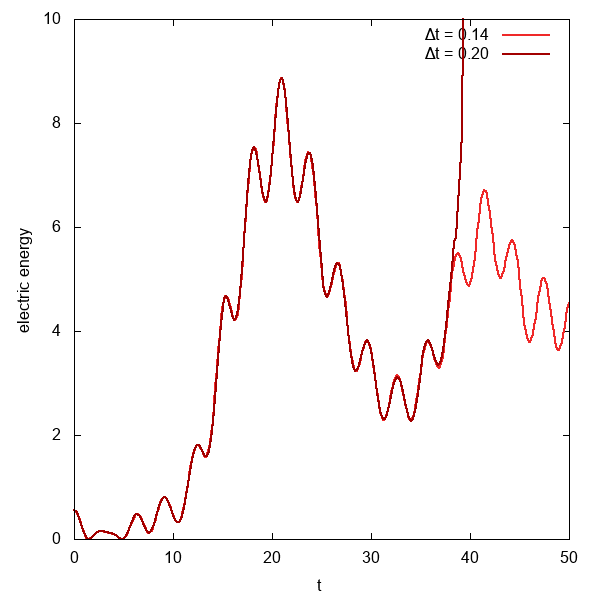} & \hspace{-0.2cm}\includegraphics[scale=0.3]{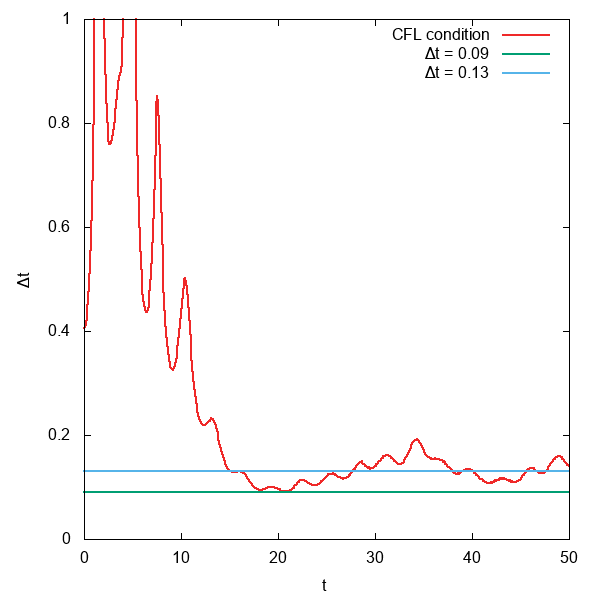} 
   \end{tabular}
\caption{Illustration of the accuracy of the CFL estimate obtained from the linear theory. History of electric energy with Lawson($RK(4,4)$) + WENO5 (left),  Lawson($RK(4,4)$) + centered scheme (middle) and history of CFL condition for Lawson($RK(4,4)$) + WENO5 case (right)}
\label{unstable}      
\end{figure}

Although being able to choose the time step size independently of the mesh is a desirable feature, it makes checking the sharpness of the CFL estimate derived in the previous section more difficult. To accomplish this, we consider the same parameters as before, except for the phase space mesh which now uses $N_x=81$ and $N_v=512$ grid points.  Then the maximum time step becomes $\Delta t=\min_n C\Delta v/\|E^n\|_{L^\infty} \approx 0.052C$ (since $\max_n \|E^n\|_{L^\infty}\approx 0.6$ and $\Delta v=16/512=0.03125$). We consider two different time steps:  $\Delta t=0.052C$ (which satisfies the linearized CFL condition) and $\Delta t = 1.4 \times 0.052C$ (which violates the linearized CFL condition). The results are shown in Figure \ref{unstable}. There the Lawson($RK(4,4)$) method has been chosen for the time discretization whereas WENO5 and centered scheme are both considered for the velocity discretization. More specifically,
\begin{itemize}
    \item for WENO5 we use $C=1.73$ (obtained from the linearized version LW5) and we compare the results obtained 
    with $\Delta t=0.09$ (satisfies the CFL condition) and $\Delta t=0.13$ (does not satisfy the CFL condition). 
    \item for the centered scheme, we use $C=2\sqrt{2}$ and we compare the results obtained with $\Delta t=0.14$ (satisfies the CFL condition) with $\Delta t=0.2$ (does not satisfy the CFL condition). 
\end{itemize}
In Figure \ref{unstable}, the time evolution of the electric energy $\|E(t)\|^2_{L^2}$ is displayed for these two velocity discretizations.  
One can observe for the time step size that satisfies the CFL condition the simulation is stable and gives the expected results, whereas for the choice that violates the CFL condition the simulation blows up. Thus, the results confirm that the CFL condition obtained by the linear theory yields a good prediction for the nonlinear Vlasov--Poisson equation. 
On the right part of Figure \ref{unstable}, the time history of the quantity $C \Delta v/ \|E^n\|_{L^\infty}$ is shown (red) together with the time step size considered for the WENO velocity discretization. The choice $\Delta t=0.13$ (blue line) is larger than the allowed time step size around $t\approx 20$, which explains the numerical instability observed at that point in time. 

\section{Numerical simulation: drift-kinetic equations \label{sec:dk}}

In this section we will consider a model motivated by the simulation of strongly magnetized plasmas, such as those found in tokamaks. In this case the dynamics is governed by gyrokinetic equations. Gyrokinetics averages out the fast oscillatory motion of the charged particles around the magnetic field lines. In a simplified slab geometry, gyrokinetic models reduce to the drift-kinetic equation. In this case the unknown $f$ depends on three cylindrical spatial coordinates $(r,\theta,z)$ and one velocity direction $v$. This model is composed of a guiding-center dynamics in the plane orthogonal to the magnetic field lines and of a Vlasov type dynamics in the direction parallel (to the magnetic field lines). In addition to its relevance in physics, it is also a good test case for stressing exponential methods. The latter is due to the fact that after some time the nonlinearity can become strong enough such that the time step size is is dictated by stability constraints (especially for high order methods).

Our goal in this section is to find a numerical approximation of $f=f(t,r,\theta,z,v)$ satisfying the following $4D$ slab drift-kinetic equation (see \cite{grandgirard})
\begin{equation}
\label{dk}
\partial_tf-\frac{\partial_\theta \phi}{r}\partial_rf+\frac{\partial_r \phi}{r}\partial_\theta f
+v\partial_zf-\partial_z\phi\partial_{v}f=0,
\end{equation}
for $(r,\theta,z,v)\in \Omega\times[0,L]\times \mathbb{R}$, $\Omega=[r_{\rm min},r_{\rm max}]\times [0, 2\pi]$.
The self-consistent potential $\phi=\phi(r,\theta,z)$ is determined by solving the quasi neutrality equation
\begin{align}
    -\left[\partial_r^2\phi+\left(\frac{1}{r}+\frac{\partial_r n_0(r)}{n_0(r)}\right)\partial_r\phi \right.&+\left.\frac{1}{r^2}\partial_\theta^2\phi\right]+\frac{1}{T_e(r)}(\phi-\langle\phi\rangle)\nonumber\\
\label{qn}
&\hspace{-3cm}=\frac{1}{n_0(r)}\int_{\mathbb{R}} fdv-1,
\end{align}
where $\langle\phi\rangle = \frac{1}{L}\int_0^L \phi(r,\theta,z)\,dz$ and the functions $n_0$ and $T_e$ depend only on $r$ and are given analytically.

In many situations the $v \partial_z f$ term yields the most restrictive CFL condition. In this setting exponential methods can be very successful as they remove the most stringent CFL condition, while still treating the remaining terms explicitly (which computationally is relatively cheap). The $\varphi$ functions can be computed in Fourier space (as has been discussed in some detail for the Vlasov--Poisson system in the previous section) or using a semi-Lagrangian approach. Exponential integrators for the drift-kinetic model have been proposed in \cite{cep}. They compare favorably to splitting schemes and have the advantage that they can be more easily adapted to different models. In \cite{cep} only a second order exponential integrator and the fourth order Cox--Matthews scheme have been considered. Due to the investigations in the present paper we now understand that this is not an ideal choice. Thus, the purpose of this section is to demonstrate that Lawson methods can be more efficient and to further corroborate the results obtained in the previous sections. The difference in stability for Lawson schemes and exponential integrators will be very evident in the numerical simulations that are presented.

\subsection{Numerical discretization}

First, we remark that $z$ is a periodic variable which motivates us to consider the Fourier transform in this direction. The corresponding frequencies are denoted by $k$. Equation \eqref{dk} then becomes
$$
\partial_t \hat{f}_k -\partial_r \widehat{\left(\frac{\partial_\theta \phi}{r}f\right)}_k+\partial_\theta \widehat{\left(\frac{\partial_r \phi}{r} f\right)}_k +vik \hat{f}_k-\partial_v\widehat{\left( \partial_z\phi \, f\right)}_k=0, 
$$
Setting $F(t, f)= \partial_r \widehat{\left(\frac{\partial_\theta \phi}{r}f\right)}-\partial_\theta \widehat{\left(\frac{\partial_r \phi}{r} f\right)} +\partial_v\widehat{\left( \partial_z\phi \,f\right)}$, this equation can be written as
$$
\partial_t \hat{f} = - vik \hat{f} + F(t, f). 
$$
This is now precisely in the form to which we can apply an exponential method. In addition, computing the required matrix functions is very efficient as all the frequencies decouple (see the corresponding discussion in section \ref{sec:vp}).

To complete the numerical scheme, one has to detail the phase space approximation. As in \cite{cep} we will use Arakawa's method to approximate the derivatives needed to compute $F$. Arakawa's method is a centered difference scheme that conserves three invariants. More details can be found in \cite{cep}.

\subsection{Numerical results \label{subsec:driftkinetic-results}}

In this section, we detail the physical parameters of the considered test case. The set up is identical to \cite{cep} (see also 
\cite{BC2013, vlasovia}). The initial value is given by 
\begin{align*}
f(t=0,r,\theta,z,v) &=
f_{\rm eq}(r,v)\left[1+\epsilon \exp\left(-\frac{(r-r_p)^2}{\delta r}\right)\cos\left(\frac{2\pi n}{L}z+m\theta\right)\right],
\end{align*}
where the equilibrium distribution is given by
\begin{equation} \label{eq:equilibrium}
f_{\rm eq}(r,v)=\frac{n_0(r)\exp(-\frac{v^2}{2T_i(r)})}{(2\pi T_i(r))^{1/2}}.
\end{equation}
The radial profiles $T_i$, $T_e$, and $n_0$ have the analytic expressions
$$
\mathcal{P}(r) = C_\mathcal{P}\exp\left(-\kappa_\mathcal{P}\delta r_{\mathcal{P}}\tanh(\frac{r-r_p}{\delta r_{\mathcal{P}}})\right), \; \mathcal{P}\in \{T_i,T_e,n_0\}
$$
with the constants defined as follows
$$
\ C_{T_i}=C_{T_e}=1,\ C_{n_0}=\frac{r_{\rm max}-r_{\rm min}}{
\int_{r_{\rm min}}^{r_{\rm max}}\exp(-\kappa_{n_0}\delta r_{n_0}\tanh(\frac{r-r_p}{\delta r_{n_0}}))\,dr}.
$$
Finally, we consider the parameters of \cite{BC2013} (MEDIUM case)
\begin{eqnarray*}
&&r_{\rm min} = 0.1,\ r_{\rm max} = 14.5,\\
&& \kappa_{n_0}= 0.055,\ \kappa_{T_i}=\kappa_{T_e}= 0.27586,\\
&&\delta r_{T_i}=\delta r_{T_e}=\frac{\delta r_{n_0}}{2}= 1.45,\ \epsilon=10^{-6},\ n=1,\ m=5,\\
&&L=1506.759067,\ r_p = \frac{r_{\rm min}+r_{\rm max}}{2},\delta r = \frac{4 \delta r_{n_0}}{\delta r_{T_i}}.
\end{eqnarray*}
and use a $v$-range of $v \in [-7.32,7.32]$.

We consider two configurations. A direct formulation, where the boundary conditions are given by
$$
f(r_{\rm min},\theta ,z,v)=f_{eq}(r_{\rm min},v) \qquad
f(r_{\rm max},\theta ,z,v)=f_{eq}(r_{\rm max},v).
$$
Note that these are not homogeneous Dirichlet boundary conditions. It is well known (and supported by \cite{cep}) 
that the Arakawa scheme works better for homogeneous boundary conditions. In addition to the direct formulation, we therefore also introduce a so-called perturbation formulation (see also \cite{vlasovia, latu2}).
First, we note that the equilibrium function $f_{\rm eq}$ defined in (\ref{eq:equilibrium}) is a steady state for our problem. We therefore divide $f$ into
$$
f(t,r,\theta ,v)=f_{eq}(r,v)+\delta f(t,r,\theta ,v).
$$
With this formulation, our problem (\ref{dk}) becomes
$$
\partial_t\delta f+\frac{E_\theta}{r}\partial_r (f_{eq} + \delta f)-\frac{E_r}{r}\partial_\theta \delta f+v\partial_z\delta f+E_z \partial_v (f_{eq} + \delta f)=0,
$$
where $E_\theta=-\partial_\theta \phi$, $E_r=-\partial_r \phi$ and $E_z=-\partial_z \phi$. Expanding the various terms we obtain
$$
\partial_t\delta f+\frac{E_\theta}{r}\partial_r\delta f-\frac{E_r}{r}\partial_\theta \delta f+v\partial_z\delta f+ E_z \partial_v \delta f
+\frac{E_\theta}{r}\partial_r f_{eq} + E_z \partial _v f_{eq}=0
$$
which can be written as
$$
\partial_t\delta f + v\partial_z \delta f -F(\delta f) +\frac{E_\theta}{r}\partial_r f_{eq} +E_z\partial _v f_{eq} = 0.
$$
Note that the equation is very similar to equation (\ref{dk}). We, however, have obtained two additional source terms, which depend on the equilibrium distribution $f_{eq}$ as well as on the electric field.
Furthermore, the right hand side of the quasi-neutrality equation (\ref{qn}) becomes
$$
\frac{1}{n_0}\int f_{eq} \,dv+\frac{1}{n_0}\int \delta f \,dv-1 = \frac{1}{n_0}\int \delta f \,dv.
$$
Due to the similarity of the direct formulation and the perturbation formulation, the same code can 
be used for both by simply exchanging the right hand side of the quasi-neutrality equation, changing the boundary conditions, and adding the appropriate source terms. Thus, to implement the exponential integrator we consider the following equation 
$$ \partial_t\delta f + v\partial_z \delta f  = F_{pert}(\delta f),  $$
with
$$ F_{pert} (\delta f)= F(\delta f) -\frac{E_\theta}{r}\partial_r f_{eq} - E_z\partial _v f_{eq},  $$
and proceed as before (with $F$ replaced by $F_{pert}$). The space discretization of the source terms can be done either analytically or using a numerical approximation. In our implementation we have used standard centered differences. The Arakawa scheme that is used to discretize $F(\delta f)$ now employs homogeneous Dirichlet boundary conditions for $\delta f$ in the $r$-direction.

We have seen in section \ref{sec:vp} that for the Vlasov--Poisson equation we can derive a constraint on the time step size which ensures stability. 
For Lawson methods this also gives a good estimate in practice. However, for exponential integrators the situation is far more complicated, see the discussion in section \ref{ode}. Thus, a natural question that arises is how large time steps can we take in practice. To do that we will employ an adaptive step size controller that uses Richardson extrapolation to obtain an error estimate. By denoting a time step as follows $f^{n+1} = \varphi_{\Delta t_n}(f^n)$ and considering 
$\tilde{f}^{n+1} = \varphi_{\Delta t_n/2}\circ \varphi_{\Delta t_n/2}(f^n)$
we can construct the Richardson extrapolated numerical solution of a method of order $p$ as follows
$f_R^{n+1} = (2^{p+1} \tilde{f}^{n+1} - f^{n+1})/(2^{p+1} -1)$, which turns out 
to be an approximation of order $(p+1)$ of the exact solution.  
Then, it is possible to determine an estimate of the local error 
$e_{n+1}$ 
of the time integrator through the following expression 
$$
e_{n+1} = \left\Vert f_R^{n+1} - f^{n+1}\right\Vert_{L^{\infty}} + {\cal O}(\Delta t_n^{p+2}), 
$$
where the $L^\infty$ norm is considered in the $r, \theta, z, v$ variables. If the estimate for the error $e_{n+1}$ is larger than a specified $\text{tol}$ we reject the step and start again from time $t_n$. Otherwise, the step is accepted and we proceed with the time integration. In either case we then determine the new step size $\Delta t_{new}$ such that the local error is smaller than the tolerance. That is, we choose 
\begin{equation} 
\label{compute_dt}
\Delta t_{new}=s \Delta t_{n}\left(\frac{\text{tol}}{e_{n+1}}\right)^{1/(p+1)}, 
\end{equation}
where $\text{tol}$ is the prescribed tolerance, $p$ is the order of the method, and $s=0.8$ is a safety factor. This process is very well established in the literature and we refer the interested reader to \cite{gustafsson1988,gustafsson1994,soderlind2002,soderlind2006,lukas}. 
Other strategies can also be considered such as embedded Lawson or exponential methods (see \cite{gni, balac}).  
Such methods may be more efficient but we restrict ourselves here to the strategy based on Richardson extrapolation 
since it can be applied to any time integrator.

An interesting property of this adaptive step size controller is that it forces the time step size to satisfy the stability constraint of the numerical method. This is perhaps surprising at first sight since the scheme only controls the local error. However, numerical instability are characterized by error amplification as integration proceeds in time. Thus, a single step can violate the stability constraint, but later on the error amplification increases the local error in such a way that the adaptive step size controller is forced to reduce the time step size. Thus, the controller ensures that we obtain a stable numerical simulation for which the local error is below the specified tolerance.

This procedure allows us to perform a fair comparison between Lawson methods and exponential integrators. Since we are mainly interested in the stability of the methods and, particularly in the nonlinear regime, prescribing a stringent tolerance is infeasible in any case, we will choose a relatively large tolerance for our simulation ($\text{tol}=10^{-2}$ for the perturbation formulation). To avoid the problem of too large time steps at the beginning of the simulation, where accuracy and not stability dictates the time step, we limit the maximal step size to 
$\Delta t=11$ (coarse) and $\Delta t=10$ (fine) for second order methods, $\Delta t=30$ for third order methods and Lawson($RK(3, 2) \; best$), 
and $\Delta t=40$ for fourth order methods.

\begin{figure}[h]
	\centering
	\begin{subfigure}[b]{0.48\textwidth}
        \centering \includegraphics[width=\textwidth]{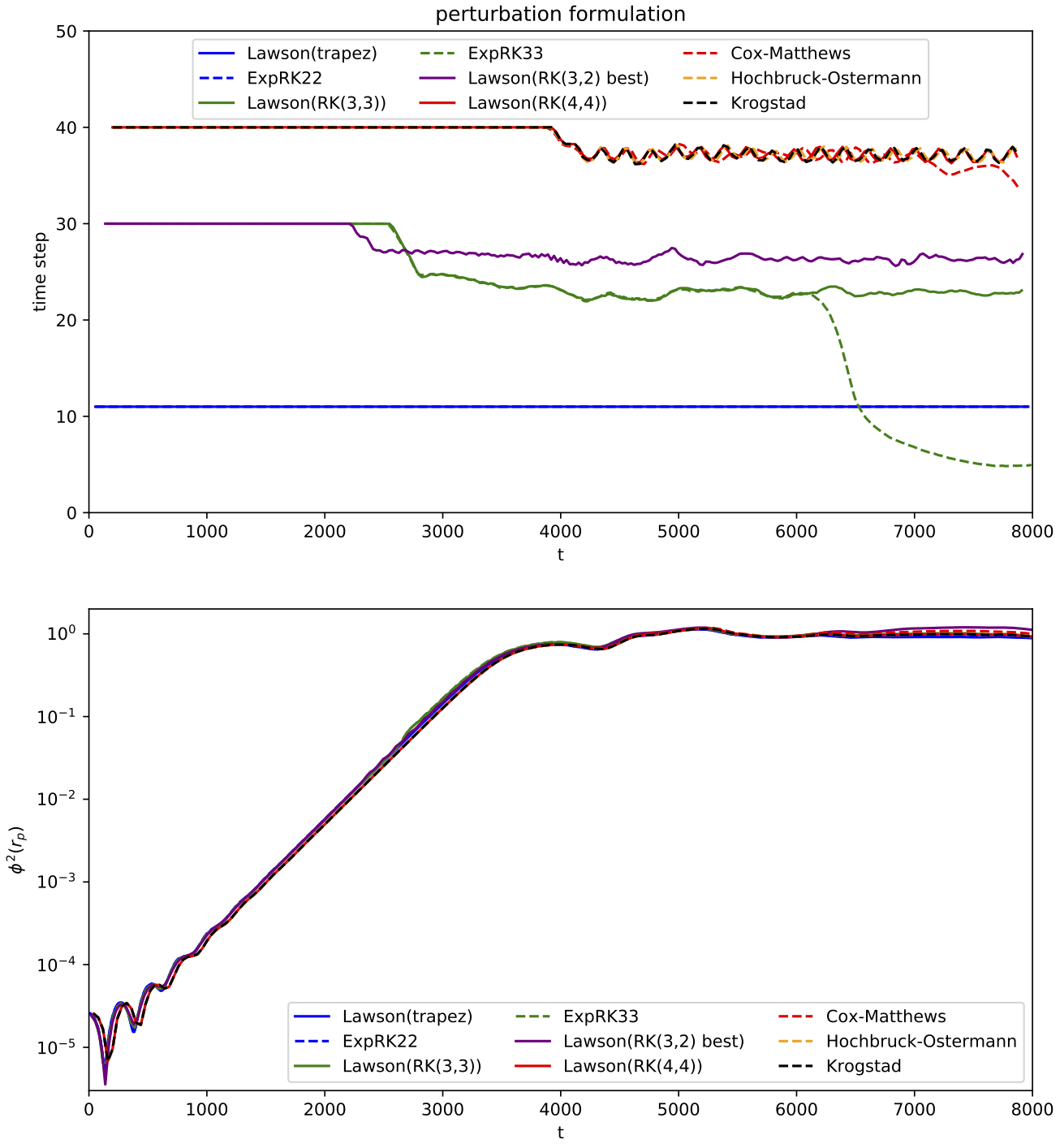}
	\end{subfigure}
	\begin{subfigure}[b]{0.48\textwidth}
        \centering \includegraphics[width=\textwidth]{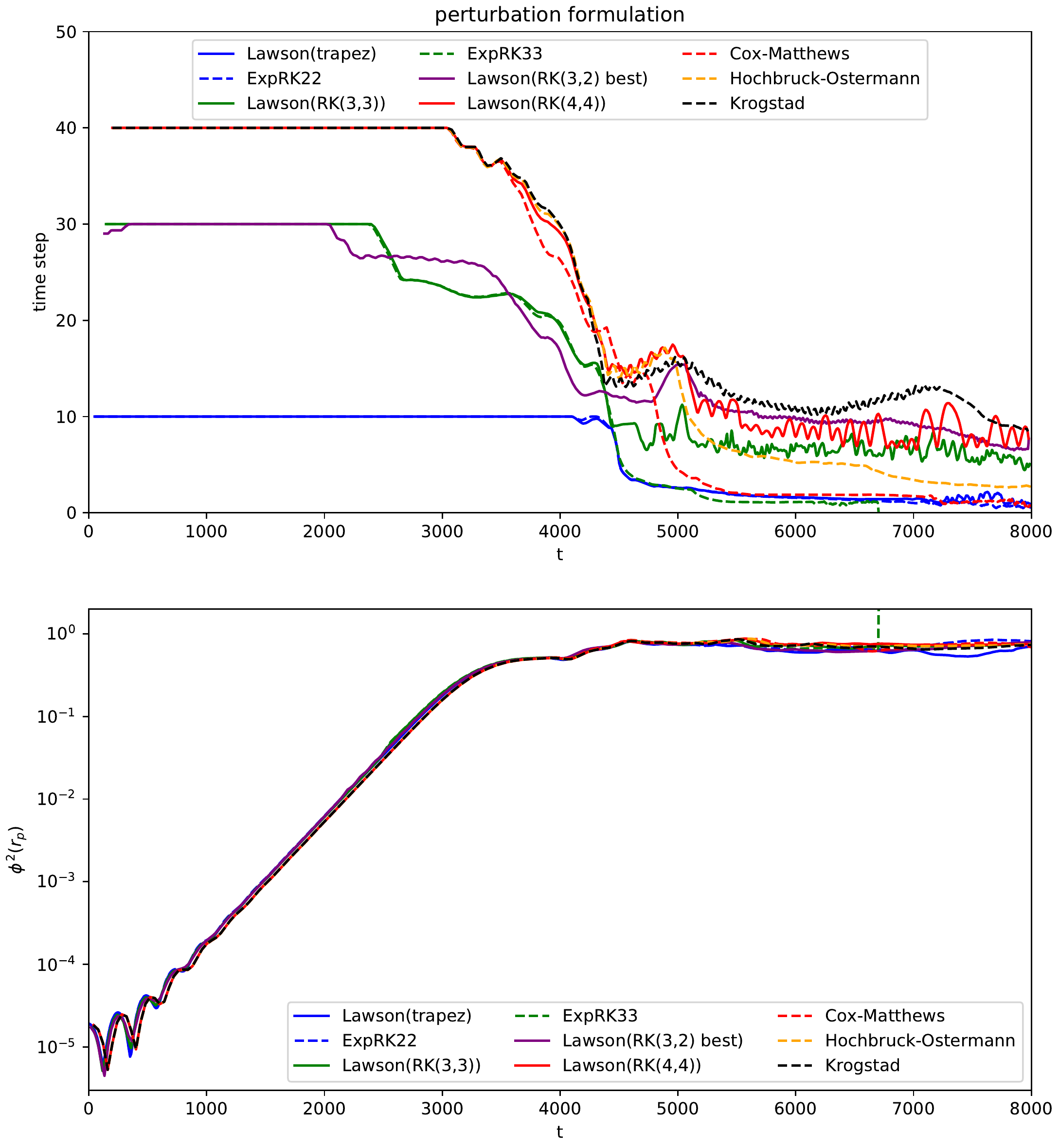}
	\end{subfigure}
    \caption{Numerical simulation for a number of Lawson methods and exponential integrators for the drift-kinetic model (perturbation formulation). The upper plots show the time step size as a function of time. The lower plots show the time evolution of the electric energy. The configuration on the left uses $32 \times 32 \times 32 \times 64$ grid points and the configuration on the right uses $64 \times 64 \times 64 \times 128$ grid points.}\label{fig:driftkinetic-pert1}
\end{figure}

To evaluate the performances of the different time integrators, we consider the time evolution of the electric energy defined by
$$
{\cal E}(t) = \left( \int_{0}^{L} \int_0^{2\pi} \phi^2(t, r_p, \theta, z ) \,d\theta dz \right)^{1/2}, \;\;\; \mbox{ with } r_p=\frac{r_{\min}+r_{\max}}{2}, 
$$
as well as the time evolution of the total mass and the total energy 
\begin{eqnarray*}
{\cal M}(t) &=& \int_{r_{\min}}^{r_{\max}} \int_{0}^{L} \int_0^{2\pi} \int_{\mathbb R} f(t, r, \theta, z, v ) \,dv d\theta dz dr, \nonumber\\ 
{\cal N}(t) &=&\int_{r_{\min}}^{r_{\max}} \int_{0}^{L} \int_0^{2\pi} \int_{\mathbb R} \frac{v^2}{2} f(t, r, \theta, z, v ) \,dv d\theta dz dr \nonumber\\ 
&& + \int_{r_{\min}}^{r_{\max}} \int_{0}^{L} \int_0^{2\pi} \int_{\mathbb R} f(t, r, \theta, z, v )\phi(t, r, \theta, z) \,dv d\theta dz dr. 
\end{eqnarray*}
The numerical results for the perturbation formulation are given in Figure \ref{fig:driftkinetic-pert1}. There
the time history of the electric energy and the time step size as of function of time for two different 
discretizations in phase space, $32 \times 32 \times 32 \times 64$ and  $64 \times 64 \times 64 \times 128$ grid points, are shown. 
We first see that all the time integrators agree very well; that is, we see an initial exponential growth (the rate is in good agreement with the linear theory; see, for example, \cite{BC2013}) 
in the electric energy. This phase is followed by saturation at very similar levels for all numerical methods used. 
We also observe that all exponential integrators, except the method of Krogstad, are forced to reduce their step size after time $t\approx 5000$ for the fine, \ie $64\times 64 \times 64 \times 128$ grid points, case. 
This is particularly drastic for ExpRK33 and the Cox--Matthews method which suffer from stability issues even for the coarse discretization. In general, for the finer space discretization (see the right plot in Figure \ref{fig:driftkinetic-pert1}), the problem becomes significantly more severe. It is also worth mentioning that ExpRK33 leads to unstable results in spite of the step size controller, which clearly highlights the unstable nature of that integrator. 
Neither of the Lawson schemes have similar issues and the size of the stability domain on the imaginary axis gives a good indication of the relative time steps this methods can take. We also note that at later times Lawson schemes are able to take significantly larger time steps compared to exponential integrators. 
The only exponential integrator that performs well in this regime is the method of Krogstad. Thus, the numerical results agree well with what we would expect based on the theoretical analysis.

\begin{figure}[h]
	\centering
	\begin{subfigure}[b]{0.48\textwidth}
        \centering \includegraphics[width=\textwidth]{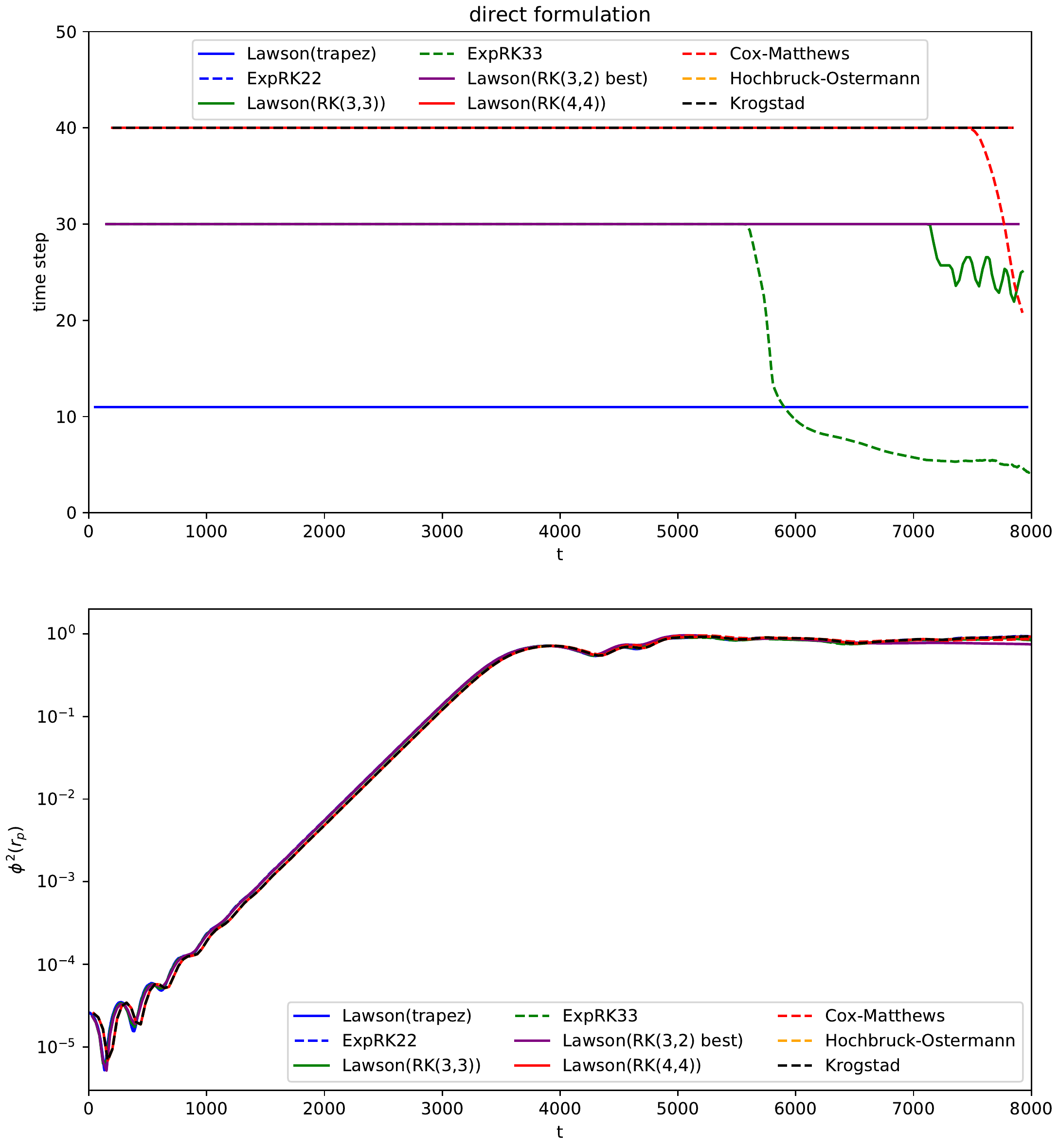}
	\end{subfigure}
	\begin{subfigure}[b]{0.48\textwidth}
                \centering \includegraphics[width=\textwidth]{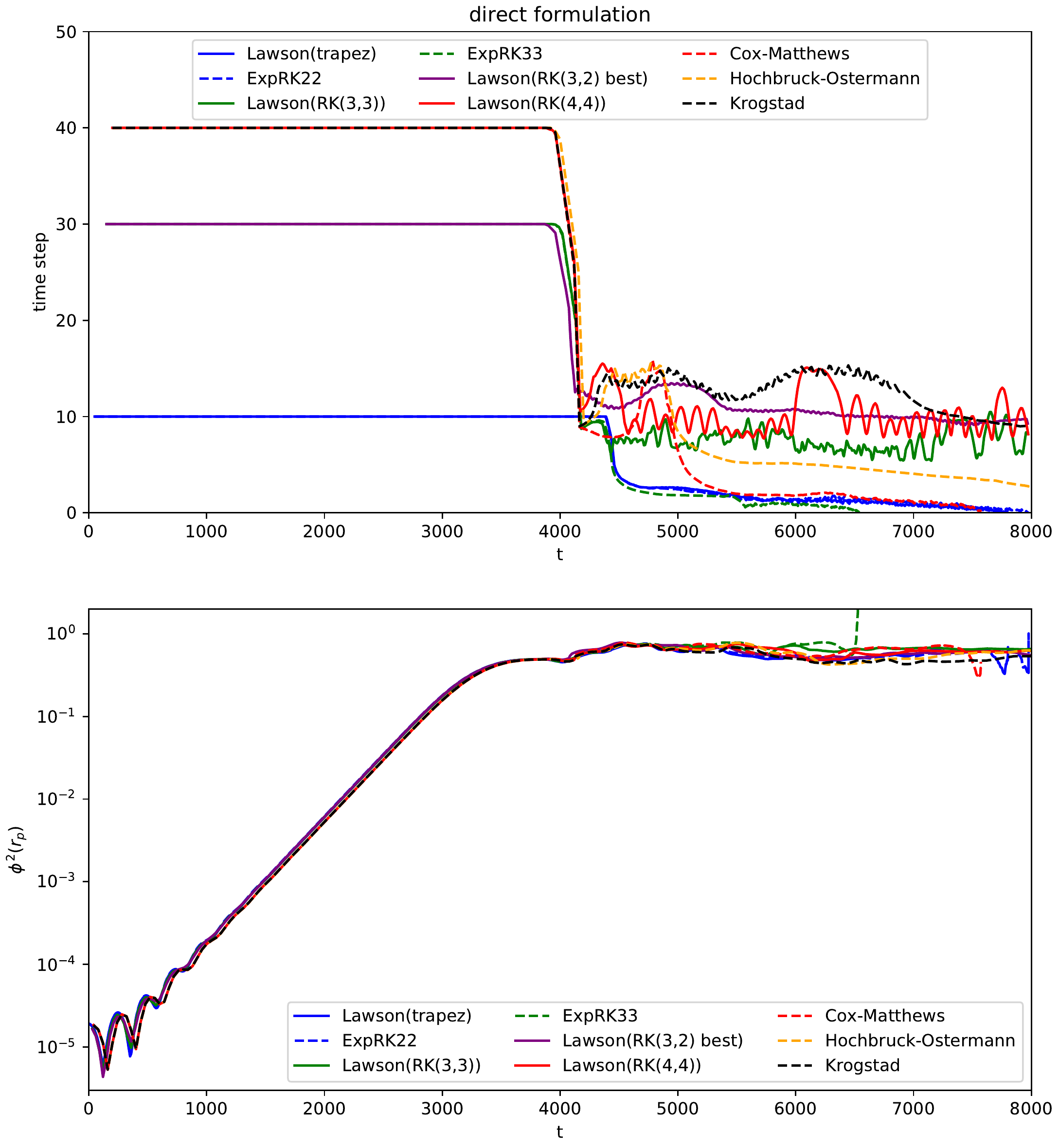}
	\end{subfigure}
	\caption{Numerical simulation for a number of Lawson methods and exponential integrators for the drift-kinetic model (direct formulation). The upper plots show the time step size as a function of time. The lower plots show the time evolution of the electric energy. The configuration on the left uses $32 \times 32 \times 32 \times 64$ grid points and the configuration on the right uses $64 \times 64 \times 64 \times 128$ grid points.}\label{fig:driftkinetic-pert0}
\end{figure}

The corresponding numerical results using the direct formulation are shown in Figure \ref{fig:driftkinetic-pert0}. The situation for the direct formulation is very similar to the perturbation formulation, even if one 
can observe that the time steps are slightly larger than in the perturbation formulation. 
One explanation comes from the fact that the relative error computed in the perturbation 
case involves the norm of $\delta f$ which can be quite different from the norm 
of $f$ so that equation \eqref{compute_dt} leads to different value of the time step 
even if the accuracy of the solution is the same.

\begin{figure}[h]
	\centering
        \centering \includegraphics[width=\textwidth]{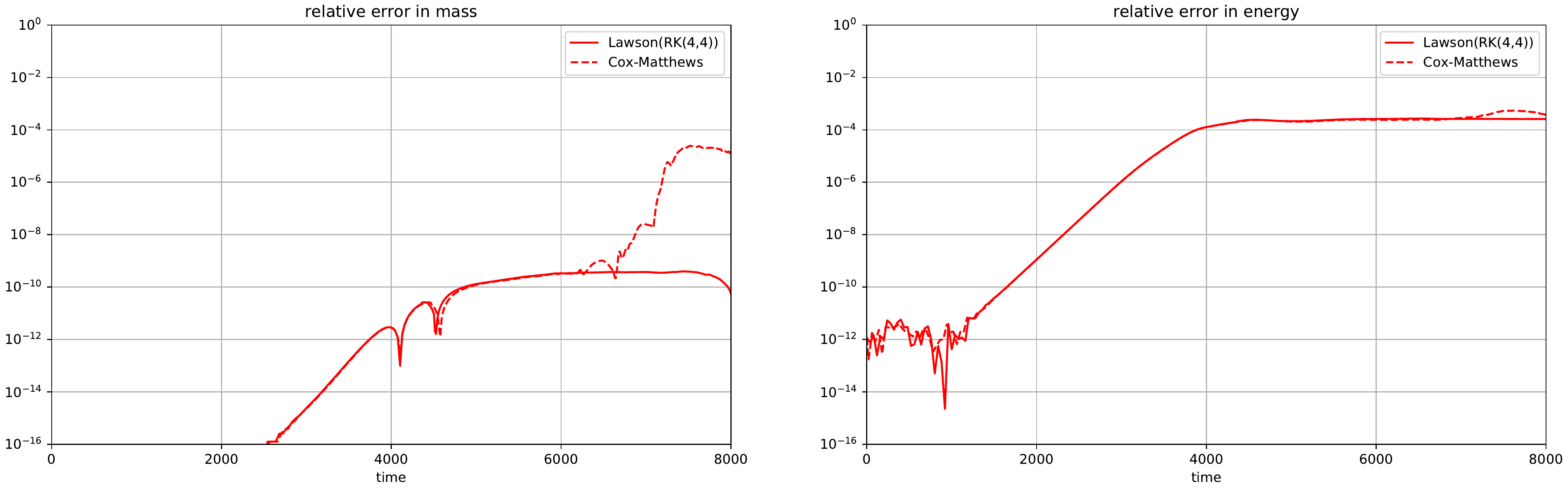}
    \caption{Numerical simulation for Lawson(RK(4,4)) and the Cox--Matthews method for the drift-kinetic model (perturbation formulation). Left: time history of the error in total mass. 
	Right: time history of the error in total energy. The configuration uses  $64 \times 64 \times 64 \times 128$ grid points.}
	\label{fig:mass_energy}
\end{figure}

In Figure \ref{fig:mass_energy}, the time history of the relative error of the total mass and of the total energy 
are displayed with the phase space discretization $64\times 64\times 64\times 128$ and 
using Lawson($RK(4,4)$) and Cox--Matthews time integrators (the perturbation formulation is used here). Since mass is a linear invariant it is preserved, up to machine precision, by the exponential integrator, see \cite{le2015KP}. We also observe good conservation of energy even in the nonlinear phase, which confirms the excellent behavior of the methods. 

\begin{figure}[h]
	\centering
    \includegraphics[width=\textwidth]{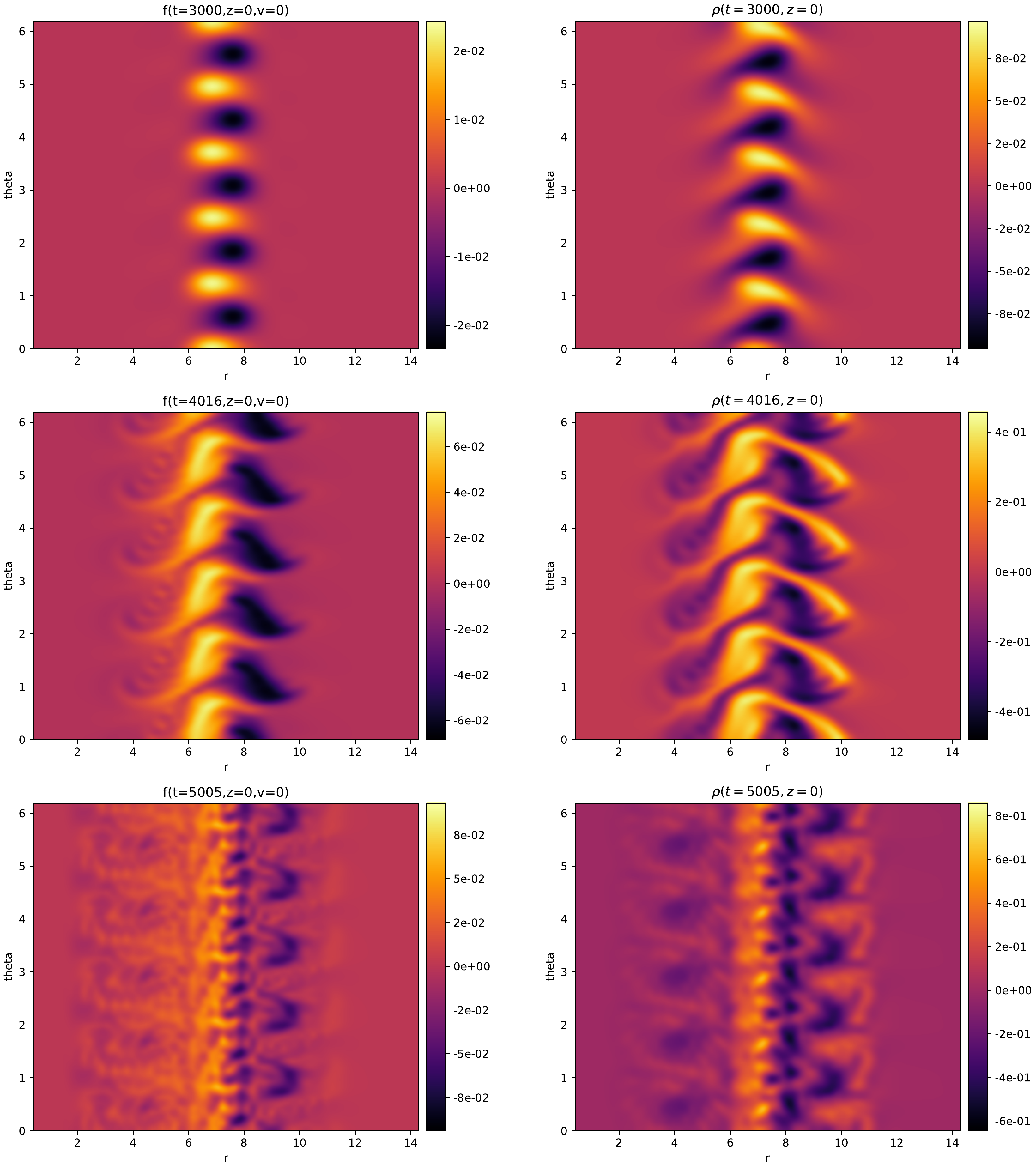}
    \caption{A slices at $(z,v)=(0,0)$ of the distribution function  (on the left) and a slice at $z=0$ of the density (on the right) are shown for times $t=3000$, $4000$, and $5000$. The Lawson($RK(4,4)$) scheme, in the configuration described in section \ref{subsec:driftkinetic-results}, with $64\times64\times64\times128$ grid points is used. \label{fig:snapshots-lrk44}}
\end{figure}

\begin{figure}[h]
	\centering
    \includegraphics[width=\textwidth]{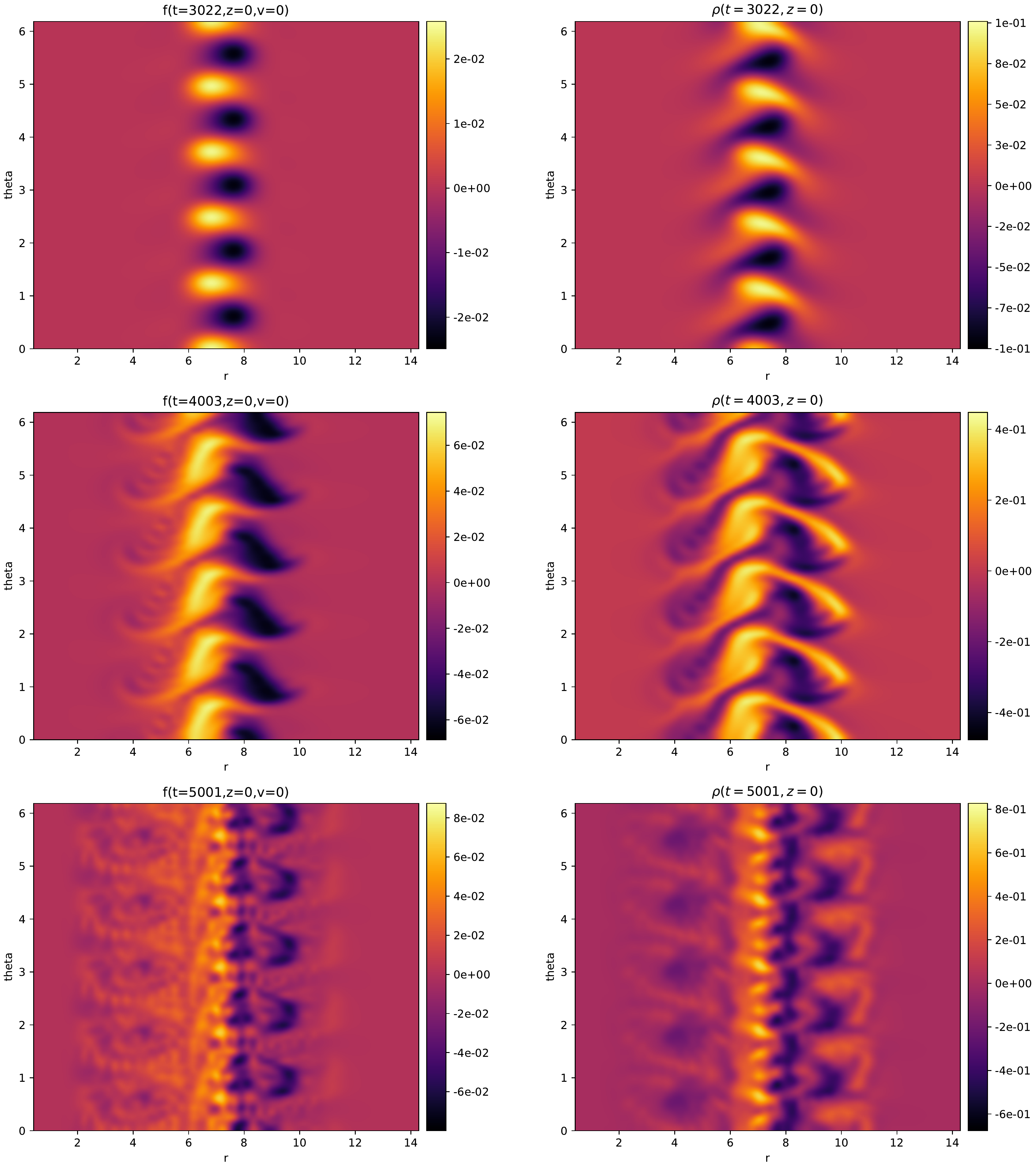}
    \caption{A slices at $(z,v)=(0,0)$ of the distribution function  (on the left) and a slice at $z=0$ of the density (on the right) are shown for times $t=3000$, $4000$, and $5000$. The Cox--Matthews scheme, in the configuration described in section \ref{subsec:driftkinetic-results}, with $64\times64\times64\times128$ grid points is used. \label{fig:snapshots-cm}}
\end{figure}

\begin{figure}[h]
	\centering
    \includegraphics[width=\textwidth]{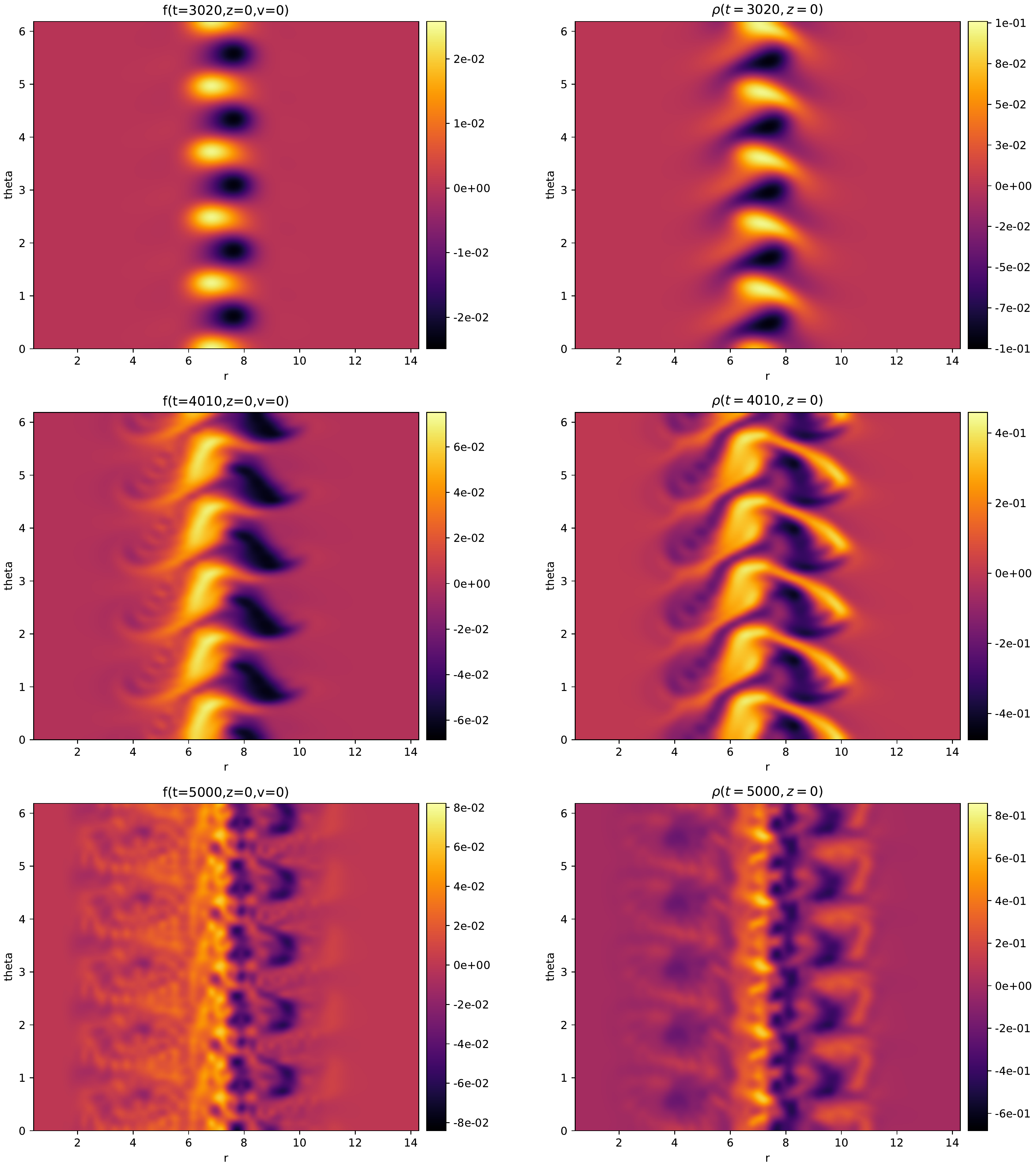}
    \caption{A slices at $(z,v)=(0,0)$ of the distribution function  (on the left) and a slice at $z=0$ of the density (on the right) are shown for times $t=3000$, $4000$, and $5000$. The Lawson($RK(4,4)$) scheme with a tolerance of $10^{-5}$ per unit step and $64\times64\times64\times128$ grid points is used. \label{fig:snapshots-ref}}
\end{figure}

Finally, we show slices of the distribution function and the density at different times. The simulation in Figure \ref{fig:snapshots-lrk44} is conducted with the Lawson($RK(4,4)$) scheme and the simulation in Figure \ref{fig:snapshots-cm} with the Cox--Matthews scheme (both 
with the perturbation formulation). In both cases the configuration of Figure \ref{fig:driftkinetic-pert1} and the fine space resolution has been employed. As comparison, a reference solution computed with the Lawson($RK(4,4)$) scheme and a step size controller that keeps the error below $10^{-5}$ per unit time step is shown in Figure \ref{fig:snapshots-ref}. We remark that all simulations show good agreement: the $m=5$ modes in the $\theta$ 
direction are recovered and after initial growth of the unstable mode,
 we can observe a shearing of the structures and the appearance of small scale structures which are typical for the nonlinear phase. 

\section*{Acknowledgement}

We would like to thank David C. Seal (U.S. Naval Academy) and Sigal Gottlieb (University of Massachusetts, Dartmouth) for the helpful discussion.

This work has been carried out within the framework of the EUROfusion Consortium and has received funding from the Euratom research and training programme 2014- 2018 and 2019-2020 under grant agreement No 633053. The views and opinions expressed herein do not necessarily reflect those of the European Commission. The work has been supported by the French Federation for Magnetic Fusion Studies (FR-FCM) and by the Austrian Science Fund (FWF): project number P 32143-N32.

\section{Appendix}

\subsection{Butcher tableaus}
\label{butcher}
In this section, we write down the different numerical methods used in this work. As in section \ref{sec:expint}, 
we consider the following equation 
$$
\dot{u} = A u + F(u), 
$$
where $A$ is a matrix and $F$ a general nonlinear function of $u$. 
The Butcher tableaus for the Lawson integrators used in the main text are stated in this section. 
A Lawson method is uniquely determined by the underlying (explicit) Runge--Kutta methods and can be written as follows  
$$
\begin{aligned}
        u^{(\ell)} &= e^{c_\ell \Delta t A}u^n + \Delta t\sum_{j=1}^s a_{\ell, j} e^{-(c_j-c_\ell)\Delta t A} F(u^{(j)}),  \\
    u^{n+1} &= e^{\Delta t A}u^n + \Delta t\sum_{j=1}^s    b_j e^{(1-c_j)\Delta t A} F(u^{(j)}), 
  \end{aligned}
$$
where the coefficients  $a_{\ell, j}$ and $b_j$ are given by the Butcher tableaus. 
The Butcher tableaus for the \textit{RK(2,2) best}, \textit{RK(3,3)} (the classic method of order 3), 
and \textit{RK(4,4)} (the classic method of order 4) are shown in Table \ref{rks}.

\begin{table}[h]
\centering
\begin{tabular}
{c|cccc}
$0$\\
$\frac{1}{2}$ & $\frac{1}{2}$\\
$\frac{1}{2}$ &$0$ &$\frac{1}{2}$ \\
$1$& $0$& $0$& $1$\\
\hline
& $\frac{1}{6}$ &$\frac{1}{3}$ &$\frac{1}{3}$ &$\frac{1}{6} $
\end{tabular}
\hspace{1cm}
\begin{tabular}
{c|ccc}
$0$\\
$\frac{1}{2}$ & $\frac{1}{2}$\\
$1 $ &$-1$ &$2$ \\
\hline
& $\frac{1}{6}$ &$\frac{2}{3}$ &$\frac{1}{6}$ 
\end{tabular}
\hspace{1cm}
\begin{tabular}
{c|ccc}
$0$\\
$\frac{1}{2}$ & $\frac{1}{2}$\\
$\frac{1}{2}$              &$0$ &$\frac{1}{2}$ \\
\hline
& $0$ &$0$ &$1$ 
\end{tabular}
    \caption{Butcher tableaus for $RK(4,4)$ (left), $RK(3,3)$ (middle) and $RK(3,2)\text{ best}$ (right).}
\label{rks}
\end{table}

A general exponential integrator can be written as 
$$
  \begin{aligned}
        u^{(\ell)} &= u^n + \Delta t\sum_{j=1}^s a_{\ell, j}(\Delta t A)\left( F(u^{(j)}) + A u^n \right) \\
    u^{n+1} &= u^n + \Delta t\sum_{j=1}^s    b_j(\Delta t A)\left( F(u^{(j)}) + A u^n \right),
  \end{aligned}
$$
where the coefficients $a_{\ell, j}(\Delta t A)$ and $b_j(\Delta t A)$ can be written as a linear combination of $\varphi_\ell$ and $\varphi_{\ell,  j}$ 
(see \cite{ei})
$$
 \varphi_\ell(z) = \frac{e^{z} - \sum_{k=0}^{\ell-1}\frac{1}{k!}z^k}{z^\ell}, 
$$ 
and we use the notations $\varphi_\ell :=\varphi_\ell(\Delta t A)$ and $\varphi_{\ell,j} := \varphi_\ell(c_j \Delta t A)$. The coefficients are collected in tableau form, see Table \ref{tab:butcher_expRK}.
\begin{table}[h]
  \centering
  \begin{tabular}{c|ccccc}
    $0$    &  & $$ & $$  & \\ 
    $c_2$    & $a_{2, 1}$ &  & $$ &  \\
    $\vdots$ & $\vdots$ & $\ddots$ & $$ &  \\
    $c_s$    & $a_{s1}$ & $\cdots$ & $a_{s, s-1}$ &  \\ \hline
    & $b_1$ & $\cdots$ & $b_{s-1}$ & $b_s$
  \end{tabular}
    \caption{Butcher tableau of a general exponential integrators}\label{tab:butcher_expRK} 
\end{table}
The Butcher tableaus for the exponential integrators used in the main text are given in Tables \ref{butcherexprk22}, \ref{butcherK}, \ref{butcherHO} and \ref{butcherCM}. 

\begin{table}[H]
  \centering
  \begin{tabular}{c|cc}
    $0$ & \\
    $1$ & $\varphi_{1,2}$ \\
    \hline
    & $\varphi_1 - \varphi_2$ & $\varphi_2$
  \end{tabular}
  \caption{Butcher tableau of ExpRK22.}
  \label{butcherexprk22}
\end{table}

\begin{table}[H]
  \centering
  \begin{tabular}{c|cccc}
    $0$           & \\
    $\frac{1}{2}$ & $\frac{1}{2}\varphi_{1,2}$ \\
    $\frac{1}{2}$ & $\frac{1}{2}\varphi_{1,3}-\varphi_{2,3}$ & $\varphi_{2,3}$ \\
    $1$           & $\varphi_{1,4}-2\varphi_{2,4}$           & $0$          & $2\varphi_{2,4}$ \\
    \hline
    & $\varphi_1-3\varphi_2+4\varphi_3$ & $2\varphi_2-4\varphi_3$ & $2\varphi_2-4\varphi_3$ & $-\varphi_2+4\varphi_3$ \\
  \end{tabular}
  \caption{Butcher tableau of the Krogstad method.}
    \label{butcherK}
\end{table}

\begin{table}[H]
  \centering
  \begin{tabular}{c|ccccc}
    $0$           & \\
    $\frac{1}{2}$ & $\frac{1}{2}\varphi_{1,2}$ \\
    $\frac{1}{2}$ & $\frac{1}{2}\varphi_{1,3}-\varphi_{2,3}$    & $\varphi_{2,3}$ \\
    $1$           & $\varphi_{1,4}-2\varphi_{2,4}$              & $\varphi_{2,4}$ & $\varphi_{2,4}$ \\
    $\frac{1}{2}$ & $\frac{1}{2}\varphi_{1,5}-2a_{5,2}-a_{5,4}$ & $a_{5,2}$       & $a_{5,2}$       & $\frac{1}{4}\varphi_{2,5} - a_{5,2}$ \\
    \hline
    & $\varphi_1-3\varphi_2+4\varphi_3$ & $0$ & $0$ & $-\varphi_2+4\varphi_3$ & $2\varphi_2-8\varphi_3$ \\
  \end{tabular}
  
    $$\begin{aligned} a_{5,2} &= \frac{1}{2}\varphi_{2,5}-\varphi_{3,4}+\frac{1}{4}\varphi_{2,4}-\frac{1}{2}\varphi_{3,5} \\ a_{5,4} &= \frac{1}{4}\varphi_{2,5}-a_{5,2} \end{aligned}
  $$
  \caption{Butcher tableau of the Hochbruck--Ostermann method.}
  \label{butcherHO}
\end{table}

\begin{table}[H]
  \centering
  \begin{tabular}{c|cccc}
    $0$           & \\
    $\frac{1}{2}$ & $\frac{1}{2}\varphi_{1,2}$ \\
    $\frac{1}{2}$ & $0$                        & $\frac{1}{2}\varphi_{1,3}$ \\
    $1$           & $\frac{1}{2}\varphi_{1,3}(\varphi_{0,3}-1)$ & $0$ & $\varphi_{1,3}$ \\
    \hline
    & $\varphi_1-3\varphi_2+4\phi_3$ & $2\varphi_2-4\varphi_3$ & $2\varphi_2-4\varphi_3$ & $4\varphi_3-\varphi_2$ \\
  \end{tabular}
  \caption{Butcher tableau of the Cox--Matthews method.}
    \label{butcherCM}
\end{table}

\subsection{WENO5 scheme}
\label{app_weno}
The different ingredients of the WENO5 scheme used in \eqref{vp_weno} are detailed here. First the fluxes are given by
$$
  \begin{aligned}
    {f}_{j+\frac{1}{2}}^+   =\ & w_0^+\left(  \frac{2}{6}f_{j-2} - \frac{7}{6}f_{j-1} + \frac{11}{6}f_{j}   \right)
                                +    w_1^+\left( -\frac{1}{6}f_{j-1} + \frac{5}{6}f_{j}   +  \frac{2}{6}f_{j+1} \right) \\
                                +  & w_2^+\left(  \frac{2}{6}f_{j}   + \frac{5}{6}f_{j+1} -  \frac{1}{6}f_{j+2} \right)
  \end{aligned}
$$
and
$$
  \begin{aligned}
    {f}_{j+\frac{1}{2}}^-   =\ & w_2^-\left( -\frac{1}{6}f_{j-1} + \frac{5}{6}f_{j}   + \frac{2}{6}f_{j+1} \right)
                                +    w_1^-\left(  \frac{2}{6}f_{j}   + \frac{5}{6}f_{j+1} - \frac{1}{6}f_{j+2} \right) \\
                                +  & w_0^-\left( \frac{11}{6}f_{j+1} - \frac{7}{6}f_{j+2} + \frac{2}{6}f_{j+3} \right).
  \end{aligned}
$$
The weights are defined through the $\beta$ coefficients 
$$
  \begin{aligned}
    \beta_0^+ &= \frac{13}{12}(\underbrace{f^+_{j-2} - 2f^+_{j-1} + f^+_{j}  }_{\Delta x^2(f''_j + \mathcal{O}(\Delta x))}))^2 + \frac{1}{4}( \underbrace{f^+_{j-2} - 4f^+_{j-1} + 3f^+_{j}}_{2\Delta  f'_j + \mathcal{O}(\Delta x^2))}  )^2 \\
    \beta_1^+ &= \frac{13}{12}( \underbrace{f^+_{j-1} - 2f^+_{j}   + f^+_{j+1}}_{\Delta x^2(f''_j + \mathcal{O}(\Delta x^2))} )^2 + \frac{1}{4}( \underbrace{f^+_{j-1} -  f^+_{j+1}}_{2\Delta x f'_j + \mathcal{O}(\Delta x^2))})^2 \\
    \beta_2^+ &= \frac{13}{12}( \underbrace{f^+_{j}   - 2f^+_{j+1} + f^+_{j+2}}_{\Delta x^2(f''_j + \mathcal{O}(\Delta x))} )^2 + \frac{1}{4}(\underbrace{3f^+_{j}   - 4f^+_{j+1} +  f^+_{j+2}}_{-2\Delta  f'_j + \mathcal{O}(\Delta x^2))})^2 \\
  \end{aligned}
$$
with 
$$
  \begin{aligned}
    \beta_0^- &= \frac{13}{12}(f^-_{j+1} - 2f^-_{j+2} + f^-_{j+3})^2 + \frac{1}{4}(3f^-_{j+1} - 4f^-_{j+2} +  f^-_{j+3})^2 \\
    \beta_1^- &= \frac{13}{12}(f^-_{j}   - 2f^-_{j+1} + f^-_{j+2})^2 + \frac{1}{4}( f^-_{j}   -  f^-_{j+2})^2 \\
    \beta_2^- &= \frac{13}{12}(f^-_{j-1} - 2f^-_{j}   + f^-_{j+1})^2 + \frac{1}{4}( f^-_{j-1} - 4f^-_{j}   + 3f^-_{j+1})^2 \\
  \end{aligned}
$$
Then, the normalized weights are 
$$
  \alpha_i^\pm = \frac{\gamma_i}{(\varepsilon + \beta_i^\pm)^2},\quad i=0,1,2, 
$$
where  $\varepsilon$ is a numerical regularization parameter set to $10^{-6}$ 
and $\gamma_0=\frac{1}{10}$, $\gamma_1=\frac{6}{10}$ 
and $\gamma_2=\frac{3}{10}$. Finally the weights are given by
$$
  w_i^\pm = \frac{\alpha_i^\pm}{\sum_m \alpha_m^\pm},\quad i=0,1,2. 
$$

\end{document}